\documentclass[12pt,twoside,english]{article}
\usepackage{amssymb,amsmath,babel}

\newtheorem{Lemma}{Lemma}[section]
\newtheorem{Theorem}[Lemma]{Theorem}

\newtheorem{Proposition}[Lemma]{Proposition}
\newtheorem{Corollaire}[Lemma]{Corollary}
\newtheorem{Definition}[Lemma]{Definition}
\newtheorem{Remarque}[Lemma]{Remark}
\newtheorem{Conjecture}[Lemma]{Conjecture}
\newtheorem{Question}[Lemma]{Question}

\newcommand{\ZZ}{\mathbb{Z}}

\newcommand{\FF}{\mathbb{F}}
\newcommand{\CC}{\mathbb{C}}

\newcommand{\GG}{\mathbb{G}}

\newcommand{\MM}{\mathbb{M}}
\newcommand{\NN}{\mathbb{N}}

\newcommand{\TT}{\mathbb{T}}

\newcommand{\rg}{\mathop{\rm rank}\nolimits}
\newcommand{\hPsi}{\boldsymbol{\widehat{{\Psi}}}}

\newcommand{\bfs}{\boldsymbol{s}}
\newcommand{\bfz}{\boldsymbol{\xi}}
\newcommand{\bfh}{\boldsymbol{h}}
\newcommand{\bfe}{\boldsymbol{E}}
\newcommand{\bfj}{\boldsymbol{J}_{\gamma}}
\newcommand{\bfl}{\boldsymbol{L}_{\gamma}}

\newcommand{\bsb}{\boldsymbol}
\newcommand{\res}{\mathbf{Res}}

\newcommand{\eg}{\emph{e.g.}}

\newcommand{\sqm}[4]{
\left(\begin{array}{ll}#1 & #2 \\ #3 & #4\end{array}\right)}
 
\newcommand\CVD{{\hfill\hfil{\lower 2 pt\hbox{\vrule\vbox to 7pt 
{\hrule width 6pt\vfill\hrule}\vrule}}}\vskip 0.5cm}

\def\carlitz{{{\text{\tiny Car}}}}

\def\tpi{{\tilde{\pi}}}

\def\GL{{\bf GL}}
    
\let\le=\leq  
\let\ge=\geq

\title{Drinfeld $A$-quasi-modular forms.}
\author{Vincent Bosser, Federico Pellarin}

\begin{document}

\maketitle

\section{Introduction}
The aim of this article is twofold: first, improve the multiplicity estimate
obtained by the second author in \cite{Pe} for Drinfeld quasi-modular forms;
and then, study the structure of certain
algebras of {\em almost-$A$-quasi-modular forms}, which already appeared in \cite{Pe}.

In order to motivate and describe more precisely our results, let us introduce
some notation. Let $q=p^e$ be a power of a prime number $p$ with $e>0$ an integer,
let $\FF_q$ be the finite field with $q$ elements. Let $\theta$
be an indeterminate over $\FF_q$, and write
$A=\FF_q[\theta]$, $K=\FF_q(\theta)$. Let $|.|$ be the absolute value on $K$ defined
by $|x|=q^{\deg_{\theta}x}$, and denote by $K_{\infty}=\FF_q((1/\theta))$ the completion
of $K$ with respect to $|.|$, by $K_{\infty}^{\text{\tiny alg}}$ an algebraic closure of $K_{\infty}$,
and  by $C$ the completion of $K_{\infty}^{\text{\tiny alg}}$
for the unique extension of $|.|$ to $K_{\infty}^{\text{\tiny alg}}$.

Let us denote by $\Omega$ the rigid analytic space $C\setminus K_{\infty}$ and by
$\Gamma:=\GL_2(A)$ the group of $2\times 2$-matrices with determinant in $\FF_q^*$, having
coefficients in $A$. The group $\Gamma$ acts on $\Omega$ by homographies.
In this setting, we can define {\em Drinfeld modular forms}
and {\em Drinfeld quasi-modular forms} for $\Gamma$ in the usual way (see \cite{BP} or
Section \ref{preliminaries} below for a definition). One of
the problems considered in this paper is to prove a {\em multiplicity estimate}
for Drinfeld quasi-modular forms, that is, an upper bound for the vanishing order at infinity
of such forms, as a function of the {\em weight} and the {\em depth}. In \cite{BP2} and \cite{Pe}, the following
conjecture is suggested
($\nu_{\infty}(f)$ denotes the vanishing order of $f$ at infinity, see Section~\ref{theproof} for the definition):

\begin{Conjecture}\label{conjecture:mult}
There exists a real number $c(q)>0$ such that, for all
non-zero quasi-modular form $f$ of weight $w$ and depth $l\geq 1$, one has
\begin{equation}\label{conjecture:bound}
\nu_{\infty}(f)\leq c(q)\,l (w-l).
\end{equation}
\end{Conjecture}

In fact, it is plausible that we can choose $c(q)=1$ in the bound (\ref{conjecture:bound}).
We refer to \cite[\S~1]{BP2} or \cite[\S~1]{Pe} for further discussion about this question.

In the classical (complex) case, the analogue of Conjecture \ref{conjecture:mult} is actually an easy exercice
using the resultant $\res_{E_2}(f,df/dz)$ in the polynomial ring $\CC[E_2,E_4,E_6]$ (here $E_{2i}$ denotes the
classical Eisenstein series of weight $2i$).

Thus, a natural idea to attack Conjecture \ref{conjecture:mult}
is to try to mimic this easy proof. However, as explained in \cite[\S~1.2]{BP2} and \cite[\S~1.1]{Pe},
if we do this we are led to use not only the first derivative of $f$ but also its {\em higher divided derivatives},
or more precisely the sequence of its {\em hyperderivatives} $D_nf$, $n\ge 0$, as defined in \cite{BP}.
But then, due to the erratic behaviour of the operators $D_n$, obstacles arise which are not easy to
overcome, and this approach appears as unfruitful to solve conjecture \ref{conjecture:bound} 
(see \cite[\S~1.2]{BP2} and \cite[\S~1]{Pe} for more details).

Another approach to prove Conjecture \ref{conjecture:mult} was carried out in \cite{BP2}. The idea was
here to use a {\em constructive} method: namely, we have constructed explicit families of {\em extremal}
Drinfeld quasi-modular forms, and, by using a resultant argument as above (the function $df/dz$
being replaced now by a suitable extremal form), we were able to get partial multiplicity estimates
in the direction of Conjecture \ref{conjecture:mult}. Unfortunately, we could not construct 
enough families of extremal forms to prove a general estimate.
Thus, also this approach to conjecture \ref{conjecture:bound} seemed unfruitful.

Recently, a new approach was introduced, this time successfully, in \cite{Pe} to get a general
multiplicity estimate (although not optimal) toward Conjecture \ref{conjecture:mult}.
The result obtained is a bound of the form
\begin{equation}\label{fedebound}
\nu_{\infty}(f)\leq c(q)\, l^2w\max\{1, \log_q{w}\},
\end{equation}
where $c(q)$ is explicit and $\log_q$ is the logarithm in base $q$. Morever, it is also proved in \cite{Pe}
that a bound like (\ref{conjecture:bound}) holds {\em if} an
extra condition of the form $w>c_0(q)l^{5/2}$ is fulfilled ($c_0(q)$ being explicit).

One of the main results of this paper is an improvement of the bound
(\ref{fedebound}), yielding Conjecture \ref{conjecture:mult} "up to a logarithm",
namely:

\begin{Theorem}\label{estimate}
There exists a real number $c(q)>0$ such that the
following holds. Let $f$ be a non zero quasi-modular form
of weight $w$ and depth $l\geq 1$.
Then
$$
\nu_{\infty}(f)\leq c(q)\, l (w-l)\max\{1,\log_q(w-l)\}.
$$
Moreover, one can take $c(q)=252q(q^2-1)$.
\end{Theorem}

The proof of this result will be given in Section~\ref{theproof}.
It consists in a refinement of the method used in \cite{Pe}. Recall that the main idea
is to introduce a new indeterminate $t$ as in Anderson's theory of $t$-motives,
and to work with certain {\em deformations} of Drinfeld quasi-modular forms (called
{\em almost $A$-quasi-modular forms} in \cite{Pe}), on which the Frobenius $\tau:x\mapsto x^q$
acts. Roughly speaking, these forms are functions $\Omega\rightarrow C[[t]]$ satisfying certain
regularity properties, as well as transformation formulas under the action of $\Gamma$
involving {\em two factors of automorphy}. The precise definitions
require quite long preliminaries: they are collected for convenience in Section~\ref{preliminaries},
which is mostly a review of facts taken from \cite{Pe}.

Section~\ref{aformes} is devoted to the problem of clarifying the structure of almost $A$-quasi-modular forms.
More precisely, let $\TT_{>0}$ denote the sub-$C$-algebra
of $C[[t]]$ consisting of series having positive convergence radius, and
let $\widetilde{\cal M}$ denote the $\TT_{>0}$-algebra of almost $A$-quasi-modular forms.
As for standard Drinfeld quasi-modular forms, almost $A$-quasi-modular forms have a depth.
Denote by $\cal M$ the sub-algebra of $\widetilde{\cal M}$ generated by forms of zero depth.
Let $E$ denote the "false" Eisenstein series of weight $2$ and type $1$ defined in \cite{Ge}.
One can define a particular almost $A$-quasi-modular form denoted by $\bfe$ (see Section~\ref{hhee}),
which is a {\em deformation} of $E$. In Section~\ref{aformes}, we obtain the following partial
description of the structure of the algebra $\widetilde{\cal M}$ (see Theorem~\ref{theoremofstructure}
for the complete statement):

\begin{Theorem}
The $\TT_{>0}$-algebra $\cal M$ has dimension $3$ and the algebra $\widetilde{\cal M}$ has
dimension $5$. Moreover, we have $\widetilde{\cal M}={\cal M}[E,\bfe]$.
\end{Theorem}
We conjecture that ${\cal M}$ is generated by three elements that can be explicitly given (see Conjecture \ref{conj}). However, we don't know how to prove this yet.

In the very last part of Section~\ref{aformes}, we define the notion of {\em $A$-modular forms}:
they generate a sub-$\TT_{>0}$-algebra of $\widetilde{\cal M}$ denoted by $\MM$. We show (Theorem~\ref{structureofaforms}) that
the algebra $\MM$ is of finite type and dimension three over $\TT_{>0}$ and we determine explicit generators.

\newpage

\section{Preliminaries}\label{preliminaries}

This section collects the preliminaries which will be needed in the next two sections.
This is essentially a review of the paper~\cite{Pe}.

\subsection{Drinfeld modular forms and quasi-modu\-lar forms.}\label{reviewoffqm}

The now classical theory of {\em Drinfeld modular forms} started with the work of Goss (see \cite{Go1})
and was improved by Gekeler (cf. \cite{Ge}).
We recall here briefly the basic definitions and properties of
Drinfeld quasi-modular forms. The reader is referred to \cite{BP} for more details and proofs.

We will use the notations of the preceding section. For $\gamma=\sqm{a}{b}{c}{d}\in\Gamma$
and $z\in\Omega$, we will denote by
$\gamma(z)= \frac{az+b}{cz+d}$ the image of the homographic action of the matrix $\gamma$ on $z$.
We will further denote by $\tau:c\mapsto c^q$ the Frobenius endomorphism, generator of
the skew polynomial ring $C[\tau]=\mathbf{End}_{\FF_q-\text{lin.}}(\GG_a(C))$.

Let $\Phi_\carlitz:A\rightarrow C[\tau]$ be the Carlitz module, defined by
$$ \Phi_{\carlitz}(\theta)=\theta\tau^0+\tau.$$ 
Let $\tpi$ be one of its fundamental periods (fixed once for all), and let
$e_{\carlitz}:C\rightarrow C$ be the associated exponential function. We have $\ker e_{\carlitz}=\tpi A$
and the function $e_{\carlitz}$ has the following entire power series expansion, for all $z\in C$:
\begin{equation}\label{ecarlitz}
e_{\carlitz}(z)=\sum_{i\geq 0}\frac{z^{q^i}}{d_i},
\end{equation}
where, borrowing classical notations,
\begin{equation}\label{di}
d_0=1,\quad d_i=[i][i-1]^q\ldots [1]^{q^{i-1}}\ \text{for}\ i\geq 1
\end{equation}
and $[i]=\theta^{q^i}-\theta$.
We define the parameter at infinity (\footnote{Note that this parameter is sometimes denoted by $t(z)$ in the literature, \eg\ in \cite{Ge} and \cite{BP}.}) $u$ by setting, for
$z\in\Omega$,
\begin{equation*}
u=u(z):=\frac{1}{e_{\carlitz}(\tpi z)}.
\end{equation*}

We will say that a function $f:\Omega \rightarrow C$ is holomorphic on $\Omega$
if it is analytic in the rigid analytic sense, and will say that it is \textit{holomorphic at infinity}
if it is $A$-periodic (that is, $f(z+a)=f(z)$ for all $z\in \Omega$ and all $a\in A$)
and if there is a real number $\epsilon>0$ such that, for all $z\in\Omega$
satisfying $|u(z)|<\epsilon$, $f(z)$ is equal to the sum of a convergent series
\begin{equation*}%\label{uexpansion}
f(z)=\sum_{n\geq 0}f_n u(z)^n,
\end{equation*}
where $f_n\in C$.
In the sequel, we will often identify such a function with a formal series in
$C[[u]]$, thus simply writing
\begin{equation*}%\label{uexpansion}
f=\sum_{n\geq 0}f_n u^n.
\end{equation*}
%We will note $\text{Hol}(\Omega)$ the ring of holomorphic functions on $\Omega$.

We can now recall the definition of a quasi-modular form for the group $\Gamma$ \cite{BP}.

\begin{Definition}{\em 
Let $w\ge 0$ be an integer and $m\in\ZZ/(q-1)\ZZ$.
A holomorphic function $f:\Omega\rightarrow C$ is a
\emph{quasi-modular form of weight $w$ and type $m$} if there exist
$A$-periodic functions $f_0,\ldots,f_l$, holomorphic in $\Omega$ and at infinity,
such that, for all $z\in 
\Omega$
and all $\gamma=\begin{pmatrix} a&b\\ c&d\end{pmatrix}\in\Gamma$:
\begin{equation*}
f(\gamma(z))=(cz+d)^w(\det\gamma)^{-m} \sum_{i=0}^l f_i(z)
\bigl(\frac{c}{cz+d}\bigr)^i.\label{QMcondition}
\end{equation*}}
\end{Definition}
In the definition above, the functions $f_0,\ldots,f_l$ are uniquely determined by $f$,
and moreover we have $f_0=f$. When $f\not=0$, the weight $w$
and the type $m$ are also uniquely determined by $f$. If $f\not=0$, we can take $f_l\not=0$ in
the definition above and $l$ is then called the \textit{depth} of $f$. The zero function
is by convention of weight $w$, type $m$ and depth $l$ for all $w$, $m$ and $l$.
A quasi-modular form of depth $0$ is by definition a \textit{Drinfeld modular form} as in \eg\ \cite{Ge}.

We will denote by $M_{w,m}$ the $C$-vector space of Drinfeld modular forms of weight $w$
and type $m$, and by $\widetilde{M}_{w,m}^{\le l}$ the $C$-vector space of Drinfeld quasi-modular forms
of weight $w$, type $m$, and depth $\le l$. We further denote by $M$ (resp. $\widetilde{M}$)
the $C$-algebra of functions $\Omega\rightarrow C$ generated by modular forms (resp. by
quasi-modular forms). By \cite[Theorem (5.13)]{Ge} and \cite[Theorem 1]{BP}, we have
\begin{equation}\label{structurefqm}
M=C[g,h]\quad \text{and}\quad \widetilde{M}=C[E,g,h],
\end{equation}
where $E$, $g$, $h$ are three algebraically independent functions defined in \cite{Ge}.
The function $g$ is modular of weight $q-1$, type $0$ (this is a kind of normalized Eisenstein series),
$h$ is modular of weight $q+1$ and type $1$ (it is very similar to a normalized Poincar\'e series), and $E$
is quasi-modular of weight $2$, type $1$ and depth $1$.

Drinfeld quasi-modular forms of weight $w$, type $m$ and depth $l$ coincide with
polynomials $f\in C[E,g,h]$ which are homogeneous of weight $w$, type $m$, and such that $\deg_Ef=l$.
Moreover, by \cite[(8.4)]{Ge} we have the following
transformation formula for $E$:
\begin{equation}\label{e:transformation}
E(\gamma(z))= \frac{(cz+d)^2}{\det\gamma}\bigl(E(z)-\frac{1}{\tilde{\pi}}\frac{c}{cz+d}\bigr)\qquad
(z\in\Omega,\gamma\in\Gamma).
\end{equation}
Finally, we have (\cite[\S~10]{Ge})
\begin{eqnarray*}
g&=&1+\cdots\, \in A[[u^{q-1}]]\cr
h&=&-u+\cdots\, \in uA[[u^{q-1}]]\cr
E&=&u+\cdots\, \in uA[[u^{q-1}]]\cr
\end{eqnarray*}
where the dots stand for terms of higher order.

\subsection{The functions $\bfs_\carlitz$, $\bfs_1$ and $\bfs_2$.}

The notion of $A$-quasi modular forms and $A$-modular forms involves in a crucial 
way three particular {\em Anderson generating functions}, which will be
denoted by $\bfs_\carlitz$, $\bfs_1$ and $\bfs_2$ as in \cite{Pe}. We recall here the definitions
and properties that will be needed later. The results quoted here are taken from \cite{Pe} but some of them are 
already implicit in \cite{anderson}.

\subsubsection{Notations and definitions concerning formal series.}\label{defseries}

%The theory of $A$-forms initiated in \cite{Pe} is based on
%Anderson's idea to introduce a new indeterminate $t$, in order to
%define some \text{deformations} of Drinfeld quasi-modular forms.
%We will need some notation and definitions concerning formal
%series which we collect here.

Let $t$ be an indeterminate. As in \cite{Pe}, for any positive real number $r>0$
we will denote by $\TT_{<r}$ the sub-$C$-algebra of $C[[t]]$
whose elements are formal series converging for $t\in C$ with $|t|<r$,
and similarly we will write $\TT_{\le r}$ for the sub-algebra of elements
of $C[[t]]$ that converge for $t\in C$ with $|t|\le r$.
%Following [papa], we will simply write $\TT$ instead of $\TT_{\le 1}$.
We will denote by $\TT_{>0}$ the sub-$C$-algebra of $C[[t]]$
of elements $f\in C[[t]]$ that have a convergence radius $>0$.

%The theory of $A$-forms involves functions $\bsb{f}:\Omega\rightarrow C[[t]]$.

If a function $\bsb{f}:\Omega\rightarrow C[[t]]$ is given, we will
denote the image of $z$ indifferently by $\bsb{f}(z)$ or $\bsb{f}(z,t)$, the latter
notation being used when we want to stress the dependence in $t$.
If the series $\bsb{f}(z,t)$ converges at $t=t_0$ for some $z\in\Omega$, we will use
both notations $\bsb{f}(z,t_0)$ or $\bsb{f}(z)|_{t=t_0}$.

Let $\bsb{f}:\Omega\rightarrow C[[t]]$ be a function. We will call
\textit{convergence radius} (or simply \textit{radius}) of $\bsb{f}$
the supremum of the real numbers $r>0$ such that
for all $z\in\Omega$ and all $t_0\in C$ with $|t_0|<r$, the 
formal series $\bsb{f}(z,t)\in C[[t]]$ converges at $t_0$.
If we denote by $r_z$ the usual convergence radius of the series
$\bsb{f}(z,t)$ (for $z\in\Omega$ fixed), the convergence radius of $\bsb{f}$
is nothing else than $\inf\{r_z\mid z\in\Omega\}$.

If $k\in\ZZ$ is an integer and $f=\sum_{i\ge 0}f_it^i$ is an element of $C[[t]]$, we define the {\em $k$-th
Anderson's twist} of $f$ by $$ f^{(k)}:=\sum_{i\ge 0}f_i^{q^k}t^i. $$
It is straightforward to check that
if $f\in C[[t]]$ has a convergence radius equal to $r$, then
$f^{(k)}$ has a convergence radius equal to $r^{q^k}$.
Similarly, if $\bsb{f}:\Omega\rightarrow C[[t]]$ has a convergence
radius equal to $r$, then the convergence radius of $\bsb{f}^{(k)}$ is $r^{q^k}$.

Most of the functions $\bsb{f}:\Omega\rightarrow C[[t]]$ that will be considered in this paper
share regularity properties that will play an important role later. We have gathered these properties
in the following definition:

\begin{Definition}
{\em Let $\bsb{f}$ be a function $\Omega\rightarrow C[[t]]$. We say that 
$\bsb{f}$ is {\em regular} if the following properties hold.

\begin{enumerate}
\item The function $\bsb{f}$ has a convergence radius $>0$.

\item There exists $\varepsilon>0$ such that, for all $t_0\in C$, $|t_0|<\varepsilon$,
the map $z\mapsto \bsb{f}(z,t_0)$ is
holomorphic on $\Omega$.

\item
For all $a\in A$, $\bsb{f}(z+a)=\bsb{f}(z)$. Moreover, there exists $c>0$ such that for all $z\in\Omega$ with
$|u(z)|<c$ and $t$ with $|t|<c$, there is a convergent expansion
\begin{equation*}\label{uexpansion}
\bsb{f}(z,t)=\sum_{n,m\geq 0}c_{n,m}t^nu^m,
\end{equation*}
where $c_{n,m}\in C$.
\end{enumerate}}
\end{Definition}\label{defiregular}

We denote by $\bsb{\mathcal{O}}$ the set of regular functions; it is a $\TT_{>0}$-algebra.
It is plain that $\bsb{\mathcal{O}}$ contains at least all the Drinfeld quasi-modular forms. We also notice that
if $\bsb{f}$ belongs to $\bsb{\mathcal{O}}$,
then $\bsb{f}^{(k)}$ belongs to $\bsb{\mathcal{O}}$ for all $k\geq 0$.

\subsubsection{Anderson generating functions.}

Let $\Lambda\subset C$ be an $A$-lattice of rank $r\ge 1$.
We will denote by $\Phi_{\Lambda}$ the Drinfeld module associated with $\Lambda$, and
by $e_{\Lambda}:C\rightarrow C$ its exponential map.
The map $e_{\Lambda}$ has a power series expansion of the form
\begin{equation*}
e_{\Lambda}(\zeta)=\sum_{i\ge 0} \alpha_i(\Lambda)\zeta^{q^i}
\end{equation*} for lattice functions $\alpha_i$,
with $\alpha_0(\Lambda)=1$.
For any $\omega\in\Lambda$ we introduce, following
Anderson \cite{anderson}, the formal series $\bfs_{\Lambda,\omega}\in C[[t]]$ defined by

\begin{equation}\label{somega}
\bfs_{\Lambda,\omega}(t)=\sum_{i\ge 0}e_{\Lambda}(\frac{\omega}{\theta^{i+1}})t^i.
\end{equation}

In \cite{Pe} the following result is proved.

\begin{Proposition}\label{propgenerale}
Let $\Lambda$ be an $A$-lattice, and let $\omega\in\Lambda\setminus\{0\}$.
\begin{enumerate}
\item The series $\bfs_{\Lambda,\omega}(t)$ lies in $\TT_{<q}$ and its convergence radius is $q$.

\item For all $t\in C$ with $|t|<q$, we have
\begin{equation}\label{slambda}
\bfs_{\Lambda,\omega}(t)=\sum_{i\ge 0}\frac{\alpha_i(\Lambda)\omega^{q^i}}{\theta^{q^i}-t}.
\end{equation}

\item The function $t\mapsto \bfs_{\Lambda,\omega}(t)$ extends to a meromorphic
function on $C$ by means of the r.h.s. of (\ref{slambda}). It has a simple pole at
$t=\theta^{q^i}$ for all $i\ge 0$, with residue $-\alpha_i(\Lambda)\omega^{q^i}$.

\item Let us write $\Phi_{\Lambda}(\theta)=\theta\tau^0+l_1\tau+\cdots+l_r\tau^r$. Then the following relation holds:
\begin{equation*}
\sum_{k=1}^{r}l_k\bfs_{\Lambda,\omega}^{(k)}(t)=(t-\theta)\bfs_{\Lambda,\omega}(t).
\end{equation*}
\end{enumerate}
\end{Proposition}

\emph{Proof.} It is easy to show that, for $\omega \not=0$, the power series (\ref{somega})
is convergent if and only if $|t|<q$. The first assertion follows from this. For the others, see
\cite[\S~4.2.2]{Bourbaki}. \CVD

\subsubsection{The function $\bfs_{\carlitz}$.}

We take here $\Lambda=\tilde{\pi}A$ and $\omega=\tilde{\pi}$ (rank $1$-case).
We set:
$$ \bfs_\carlitz(t):= \bfs_{\tilde{\pi}A,\tilde{\pi}}.$$
In this case, $\Phi_{\Lambda}$ is the Carlitz module $\Phi_\carlitz$
and $e_{\Lambda}$ is the Carlitz exponential $e_\carlitz$, so $\alpha_i(\Lambda)=1/d_i$
by (\ref{ecarlitz}).

The main properties of the series $\bfs_{\carlitz}$ are summarized
in the following proposition.

\begin{Proposition}\label{propcarlitz}
\begin{enumerate}
\item The following expansion holds, for all $|t|<q$:
\begin{equation*}
\bfs_\carlitz(t)=\sum_{i\ge 0}\frac{\tilde{\pi}^{q^i}}{d_i(\theta^{q^i}-t)}.
\end{equation*}

\item We have: \begin{equation}\label{scarpi}
\left.(t-\theta)\bfs_\carlitz(t)\right|_{t=\theta}=-\tilde{\pi}.
\end{equation}
\item The series $\bfs_\carlitz$ satisfies the following $\tau-$difference equation:
\begin{equation}\label{carlitztau}
\bfs_\carlitz^{(1)}(t)=(t-\theta)\bfs_\carlitz(t).
\end{equation}
\item The following product expansion holds, for $t\in C$:
\begin{equation}\label{scarlitz-product}
\bfs_{\carlitz}(t)=\frac{1}{\bsb{\Omega}^{(-1)}(t)}=\frac{(-\theta)^{1/(q-1)}}{\displaystyle
\prod_{i\geq 0}(1-t/\theta^{q^i})},
\end{equation}
where $(-\theta)^{1/(q-1)}$ is an appropriate $(q-1)$-th root of $-\theta$ and $\bsb{\Omega}(t)$ is
the function defined in \cite[3.1.2]{ABP}.
\end{enumerate}
\end{Proposition}

\noindent\emph{Proof.} The points 1, 2 and 3 immediately follow from Proposition~\ref{propgenerale}.
Point 4 follows from \cite[Formula (31)]{Pe}.\CVD

\subsubsection{The functions $\bfs_1$ and $\bfs_2$.}

Let us now choose $z\in\Omega$, and consider the $A$-lattice $\Lambda_z:=A+zA$, which is of rank $2$, 
with associated exponential function $e_z:=e_{\Lambda_z}$.
By \cite{Ge}, the corresponding Drinfeld module $\Phi_z$ satisfies
\begin{equation*}
\Phi_z(\theta)=\theta\tau^0+\tilde{g}(z)\tau+\tilde{\Delta}(z)\tau^2,
\end{equation*}
where
\begin{equation}\label{normalisationdelta}
\tilde{g}(z)=\tilde{\pi}^{q-1}g(z) \quad \text{and}\quad 
\tilde{\Delta}(z)=\tilde{\pi}^{q^2-1}\Delta(z)\quad \text{with}\quad \Delta=-h^{q-1}.
\end{equation}
Here the functions $g$ and $h$ are those already introduced in Section~\ref{reviewoffqm}. We define:
\begin{equation*}
\bfs_{1}(z,t):=\bfs_{\Lambda_z,z}(t)\quad \text{and}\quad 
\bfs_{2}(z,t):=\bfs_{\Lambda_z,1}(t).
\end{equation*}

We have:
\begin{Proposition}\label{props1s2}
\begin{enumerate}
\item For all $z\in\Omega$, the series $\bfs_1(z)$ and $\bfs_2(z)$ are units in $\TT_{>0}$.

\item We have: \begin{equation}\label{periods}
\left.(t-\theta)\bfs_1(z,t)\right|_{t=\theta}=-z\quad
\text{and}\quad
\left.(t-\theta)\bfs_2(z,t)\right|_{t=\theta}=-1.
\end{equation}

\item The series $\bfs_1$ and $\bfs_2$ satisfy the following linear $\tau$-difference equations
of order~$2$:
\begin{equation}\label{taudifference}
\bfs_1^{(2)} = -\frac{\tilde{g}}{\tilde{\Delta}}\bfs_1^{(1)}+\frac{t-\theta}{\widetilde{\Delta}}\bfs_1
\qquad{\text and}\qquad
\bfs_2^{(2)} = -\frac{\tilde{g}}{\tilde{\Delta}}\bfs_2^{(1)}+\frac{t-\theta}{\widetilde{\Delta}}\bfs_2.
\end{equation}

\item The function $\bfs_2$ is a regular function of convergence radius $q$.
Moreover, there exists a real number $c>0$ such that
the following expansion holds, for $|t|<q$ and $|u(z)|<c$:
\begin{equation}\label{s2:uexpansion}
\bfs_2(z,t)=\frac{1}{\widetilde{\pi}}\bfs_{\carlitz}(t)+\sum_{i\ge 1}\kappa_i(t)u^i\in C[[u]],
\end{equation}
where $\kappa_i\in\TT_{<q}$ for all $i\ge 1$.
\end{enumerate}
\end{Proposition}

\emph{Proof.} As already remarked in \cite{Pe}, the first statement follows from the fact that for all $z\in\Omega$, the constant term of the series
$\bfs_2(z)$ and $\bfs_1(z)$ are $e_{\Lambda_z}(1/\theta)$ and $e_{\Lambda_z}(z/\theta)$ respectively,
which never vanish for $z\in\Omega$. The statements 2 and 3
follow from Proposition~\ref{propgenerale}. The last assertion is proved
in \cite[proof of Proposition~5]{Pe}. \CVD

We notice that the function $\bfs_1$ is not regular, since it is not $A$-periodic by (\ref{eqvectorial}).

%\begin{Remarque}
%If we define
%\begin{equation*}
%\hPsi:=\left(\begin{array}{cc}
%\bfs_1 & \bfs_2\cr
%\bfs_1^{(1)} & \bfs_2^{(1)}\cr
%\end{array}\right),
%\end{equation*} then formulas (\ref{taudifference}) rewrite
%$$\hPsi^{(1)}=
%\left(\begin{array}{cc}
%0 & 1\cr
%(t-\theta)/\tilde{\Delta} & -\tilde{g}/\tilde{\Delta} \cr
%\end{array}\right)
%\hPsi.$$
%Thus $\hPsi$ defines a rigid analytic trivialisation of the $t$-motive associated to $\Phi_z$.
%\end{Remarque}

A remarkable feature concerning the series $\bfs_1$ and $\bfs_2$ is that their first twists are related to
the periods of second kind of the Drinfeld
module $\Phi_z$. %as defined in \cite[\S~7]{gekeler:compositio}. More precisely,
Let $F_z:C\rightarrow C$ be the unique $\FF_q$-linear function satisfying
$$ F_z(0)=0\quad \text{and}\quad e_z^q(\zeta)=F_z(\theta\zeta)-\theta F_z(\zeta).$$
Following \cite[\S~7]{gekeler:compositio}, we define $\eta_1:=F_z(z)$ and $\eta_2:=F_z(1)$
(periods of second kind of $\Phi_z$). We further define:
\begin{equation*}
\hPsi:=\left(\begin{array}{cc}
\bfs_1 & \bfs_2\cr
\bfs_1^{(1)} & \bfs_2^{(1)}\cr
\end{array}\right).
\end{equation*}
Then we have:

\begin{Proposition}
\begin{enumerate}

\item One has \begin{equation}\label{quasiperiods}
\bfs_1^{(1)}(z,\theta)=\eta_1,\qquad
\bfs_2^{(1)}(z,\theta)=\eta_2,
\end{equation}
and
\begin{equation}\label{determinant}
\left.(t-\theta)\det\hPsi(z,t)\right|_{t=\theta}=
\det\left(\begin{array}{cc}
-z & -1\cr
\eta_1 & \eta_2 \cr
\end{array}\right)=-\frac{1}{\tilde{\pi}^qh(z)}.
\end{equation}

\item For all $z\in\Omega$, we have the following equality in $\TT_{<q}$:
\begin{equation}\label{dethPsi}
\det\hPsi(z,t)=\frac{\bfs_{\carlitz}(t)}{\tilde{\pi}^{q+1}h(z)}.
\end{equation}
\end{enumerate}
\end{Proposition}

\noindent\emph{Proof.} Formulas (\ref{quasiperiods}) follow from \cite[formula (5.3)]{gekeler:compositio},
and Formula (\ref{determinant}) follows from (\ref{periods}), (\ref{quasiperiods}) and
\cite[Theorem~6.2]{gekeler:compositio}. The relation (\ref{dethPsi}) is proved in \cite{Pe}, during the proof
of Proposition~4.\CVD

The last important property of the series $\bfs_1$ and $\bfs_2$ is their behaviour
under the action of $\Gamma$.
We need to introduce new notation.
Let us denote by $\overline{A}$ the polynomial ring $\FF_q[t]$. If $a=a(\theta)\in A$, we define
$\bar{a}:=a(t)\in\overline{A}$. Similarly, we write
$$\overline{\gamma}:=
\left(\begin{array}{cc}
\overline{a} & \overline{b}\cr
\overline{c} & \overline{d} \cr
\end{array}\right)\in \mathbf{GL}_2(\overline{A})\quad\text{if}\quad
\gamma=
\left(\begin{array}{cc}
a & b\cr
c & d \cr
\end{array}\right)\in \Gamma.
$$
We further define the vectorial map:
$$\Sigma(z,t):=\binom{\bfs_1(z,t)}{\bfs_2(z,t)}.$$

\begin{Proposition}
For all $\gamma=
\left(\begin{array}{cc}
a & b\cr
c & d \cr
\end{array}\right)\in \Gamma$ and all $z\in\Omega$, we have
the following identity of series in $\TT_{<q}$:
\begin{equation}\label{eqvectorial}
\Sigma(\gamma(z),t)=(cz+d)^{-1}\overline{\gamma}\cdot\Sigma(z,t).
\end{equation}
\end{Proposition}

\noindent\emph{Proof.}
See \cite[Lemma 2]{Pe}.\CVD
%Follows from the formula (\ref{slambda}) and the fact that the functions
%$z\mapsto\alpha_i(\Lambda_z)$ are modular of weight $q^i-1$ and type $0$ for all $i$, see \cite{Pe}.\CVD

Let us now define, for $\gamma\in\Gamma$ and $z\in\Omega$,
\begin{equation}\label{defjgamma}
\bfz:=\frac{\bfs_1}{\bfs_2},\quad
J_{\gamma}(z)=J_{\gamma}=cz+d\quad\text{and}\quad \bfj(z)=\bfj=\overline{c}\bfz+\overline{d}.
\end{equation}
By Proposition \ref{props1s2},  $\bfs_2(z)$ is a unit in $\TT_{>0}$ for all $z\in\Omega$ and $\bfz(z)$ and $\bfj(z)$ are well defined
elements of $\TT_{>0}^\times$ for all $z$. We also notice that the function $\bfj:\Gamma\times\Omega\rightarrow\TT_{>0}^\times$
is a factor of automorphy for the group $\Gamma$ since (see \cite[\S~3.2]{Pe} for a proof)

\begin{equation}\label{factorofautomorphy}
\boldsymbol{J}_{\gamma\delta}(z)=\boldsymbol{J}_{\gamma}(\delta(z))\boldsymbol{J}_{\delta}(z).\end{equation}
We further notice that $$\bfz(z)|_{t=\theta}=z,\quad \bfj(z)|_{t=\theta}=J_{\gamma}(z).$$
Now, Formula (\ref{eqvectorial}) implies
\begin{equation}\label{s2modularity}
\bfs_2(\gamma(z),t)=J_{\gamma}^{-1}\bfj \bfs_2(z,t).
\end{equation}
It turns out that this formula has a generalisation for all the twists $\bfs_2^{(k)}$, $k\geq 1$.
Let us write, for
$\gamma=
\left(\begin{array}{cc}
a & b\cr
c & d \cr
\end{array}\right)\in \Gamma$ and $z\in\Omega$,
\begin{equation}\label{l:definition}
L_{\gamma}(z):=\frac{c}{cz+d}\qquad\text{and}\qquad\bfl(z):=
\frac{\overline{c}}{(\theta-t)(\overline{c}\bfs_1+\overline{d}\bfs_2)}.
\end{equation}
Then $$\bfl|_{t=\theta}=L_{\gamma}$$ by (\ref{periods}).
We note here that this definition of $\bfl$ differs from the one in \cite[\S~3.2]{Pe}.

Define further the sequence of series $\displaystyle(g_k^{\star})_{k\geq -1}$ by:
\begin{equation*}
g_{-1}^{\star}=0,\qquad g_0^{\star}=1
\end{equation*}
and
\begin{equation*}
g_k^{\star}=g^{q^{k-1}}g_{k-1}^{\star}
+(t-\theta^{q^{k-1}})\Delta^{q^{k-2}}g_{k-2}^{\star},\ k\geq 1.
\end{equation*}
We have $g_1^{\star}=g$. Moreover, the identity $g_k^{\star}(z,\theta)=g_k(z)$ holds for all $k\geq 0$,
where $g_k$ is the Eisenstein series defined in \cite[Formula (6.8)]{Ge}.

\begin{Proposition}
For all $k\geq 0$, all $\gamma\in \Gamma$ and all $z\in\Omega$, we have:
\begin{equation}\label{twistss2}
\bfs_2^{(k)}(\gamma(z))=J_\gamma^{-q^k}\bfj\left(\bfs_2^{(k)}(z)+(t-\theta)
\frac{g^{\star}_{k-1}(z)\bfs_{\carlitz}}{\widetilde{\pi}^{q^k+1}h(z)^{q^{k-1}}}\bsb{L}_{\gamma}\right).
\end{equation}
\end{Proposition}
\emph{Sketch of proof.} We give only the main steps of the proof and refer to \cite[Proof of Proposition 4]{Pe} for
the details. We first note that for $k=0$ Formula (\ref{twistss2}) is Formula (\ref{s2modularity}).
Thanks to (\ref{dethPsi}) one then shows the identity:
\begin{equation}\label{jmodularity}
\bfj^{(1)}=\bfj\bigl(1+\frac{(t-\theta)\bfs_{\carlitz}}{\widetilde{\pi}^{q+1}
h\bfs_2^{(1)}}\bfl\bigr).
\end{equation}
Applying Anderson's twist to the formula (\ref{s2modularity}) and using the identity (\ref{jmodularity})
above, we get the formula (\ref{twistss2}) for $k=1$. The proof then goes by induction on $k$, using the formula
$$
\bfs_2^{(k)}=-\frac{\widetilde{g}^{q^{k-2}}}{\widetilde{\Delta}^{q^{k-2}}}\bfs_2^{(k-1)}
+\frac{t-\theta^{q^{k-2}}}{\widetilde{\Delta}^{q^{k-2}}}\bfs_2^{(k-2)}
$$
for $k\geq 2$, which follows from (\ref{taudifference}).\CVD

\subsection{Almost-$A$-quasi-modular forms and the functions $\bfh$ and $\bfe$.}

In this section, we first define the notion of {\em almost-$A$-quasi-modular forms}.
Then we introduce and study two particular such forms, denoted by $\bfh$ and $\bfe$ in \cite{Pe}, which are
deformations of the Drinfeld modular and quasi-modular forms $h$ and $E$, respectively.
As we will see in Section~\ref{aformes}, the functions $\bfh$ and $\bfe$ will be the basic examples of
{\em $A$-quasi-modular forms}.

\subsubsection{Almost-$A$-quasi-modular forms}\label{aaaaa}

Recall (Section~\ref{defseries}) that we denote by $\bsb{\mathcal{O}}$ the $\TT_{>0}$-algebra of regular functions.
Following \cite[\S~4.2]{Pe}, we define:

\begin{Definition}
{\em Let $\bsb{f}$ be a regular function $\Omega\rightarrow C[[t]]$. 
We say that $\bsb{f}$ is an {\em almost-$A$-quasi-modular form of weight $(\mu,\nu)$,  
type $m$ and depth $\leq l$} if there exist regular functions $\bsb{f}_{i,j}\in\bsb{\mathcal{O}}$,
$0\le i+j\le l$, such that for all $\gamma\in\Gamma$
and all $z\in\Omega$,
\begin{equation}\label{almostquasimodularforms}
\bsb{f}(\gamma(z),t)=\det(\gamma)^{-m}J_\gamma^\mu\bsb{J}_\gamma^\nu\left(\sum_{i+j\leq l}\bsb{f}_{i,j}L_\gamma^i\bsb{L}_\gamma^j\right).
\end{equation}}
\end{Definition}
For $\mu,\nu\in\ZZ,m\in\ZZ/(q-1)\ZZ$, $l\in\ZZ_{\geq 0}$, we denote by $\widetilde{\mathcal{M}}_{\mu,\nu,m}^{\leq l}$ the $\TT_{>0}$-module of
almost $A$-quasi-modular forms of weight $(\mu,\nu)$, type $m$ and depth $\leq l$. 
We have $$\widetilde{\mathcal{M}}^{\leq l}_{\mu,\nu,m}\widetilde{\mathcal{M}}^{\leq l'}_{\mu',\nu',m'}\subset\widetilde{\mathcal{M}}^{\leq l+l'}_{\mu+\mu',\nu+\nu',m+m'}.$$
We also denote by $\widetilde{\mathcal{M}}$ the $\TT_{>0}$-algebra
generated by all the almost-$A$-quasi-modular forms. It was proved in \cite{Pe} that this algebra is graded
by the group $\ZZ^2\times \ZZ/(q-1)\ZZ$ and filtered by the depths. The problem of determining the structure
of $\widetilde{\mathcal{M}}$ will be considered in Section~\ref{aformes}.

It is clear that $\widetilde{M}^{\le l}_{w,m}\subset\widetilde{\mathcal{M}}^{\le l}_{w,0,m}$ for
all $w,m$. A non-trivial example of an almost-$A$-quasi-modular form is given by the function $\bfs_2$. By
(\ref{s2modularity}) we have indeed ${\bfs}_2\in\widetilde{\mathcal{M}}^{\le 0}_{-1,1,0}$.
In the following section we introduce two other examples of almost-$A$-quasi-modular forms.

\subsubsection{The functions $\bfh$ and $\bfe$.}\label{hhee}

\begin{Definition}
{\em We define the two functions $\bfh:\Omega\rightarrow \TT_{>0}$ and
$\bfe:\Omega\rightarrow \TT_{>0}$ by
\begin{equation*}\label{defhe}
\bfh=\frac{\widetilde{\pi}h\bfs_2}{\bfs_{\carlitz}}
\qquad \text{and} \qquad \bfe=\frac{\bfh^{(1)}}{\Delta}.
\end{equation*}}
\end{Definition}
We remark that these definitions make sense since $\bfs_{\carlitz}$
is a unit in $ \TT_{>0}$ and since $\Delta(z)$ does not vanish on $\Omega$.
We also note that by (\ref{carlitztau}), this definition of $\bfe$
coincides with the definition of \cite[\S~3]{Pe}.

The main properties of $\bfh$ and $\bfe$ are summarized in the following
two propositions:

\begin{Proposition}\label{proph}
The following properties hold.
\begin{enumerate}
\item
The function $\bfh$ is regular with an infinite convergence radius, and $\bfh(z,\theta)=h(z)$
for all $z\in\Omega$.

\item
The function $\bfh$ satisfies the following $\tau$-difference equation:
\begin{equation}\label{taudifferenceh}
\bfh^{(2)}=\frac{\Delta^{q-1}}{t-\theta^q}\bigl(g\bfh^{(1)}+\Delta\bfh\bigr).
\end{equation}

\item We have the following $u$-expansion, for $|t|<q$ and $|u(z)|$ sufficiently small:
\begin{equation*}\label{h:uexpansion}
\bfh(z,t)=-u+\cdots\in u\FF_q[t,\theta][[u^{q-1}]].
\end{equation*}

\item
The function $\bfh$ satisfies the following functional equations,
for all $z\in\Omega$ and all
$\gamma\in \Gamma$:
\begin{equation*}
\bfh(\gamma(z))=\frac{J_{\gamma}^q\bfj}{\det\gamma}\bfh(z).
\end{equation*}
\end{enumerate}
\end{Proposition}

\noindent\emph{Proof.} It is easily seen that the product occuring in (\ref{scarlitz-product})
expands as a series with infinite convergence radius. Since $\bfs_2$ is regular of radius $q$
by Proposition~\ref{props1s2}, it follows that the function $\bfs_2/\bfs_{\carlitz}$, hence $\bfh$, is regular of
radius $\ge q$. Moreover, we have $\bfh(z,\theta)=h(z)$ by (\ref{scarpi}) and (\ref{periods}).
The equation (\ref{taudifferenceh}) follows easily from the definition of $\bfh$ and from
the $\tau$-difference equations (\ref{carlitztau}) and (\ref{taudifference}) satisfied by $\bfs_{\carlitz}$
and $\bfs_2$ (taking the relations (\ref{normalisationdelta}) into account).
The point 3 of the proposition follows from \cite[Proposition 5]{Pe} and \cite[Lemma 16]{Pe}.
The point 4 is an immediate consequence of the formula (\ref{s2modularity}) and
of the fact that $h\in M_{q+1,1}$. It remains to prove that the convergence radius of $\bfh$ is infinite.
Rewrite the formula (\ref{taudifferenceh}) as follows:
\begin{equation}\label{otherequation}
\bfh=-\frac{g}{\Delta}\bfh^{(1)}+\frac{t-\theta^q}{\Delta^{q}}\bfh^{(2)}.
\end{equation}
Let $z$ be fixed, and let $r\geq q$ be the radius of $\bfh(z,t)$.
Then the left hand side of (\ref{otherequation}) has radius $r$, but the right hand side has obviously
a radius $r^q$ because of the twists. Hence $r=\infty$. \CVD

\begin{Proposition}\label{propE}
The following properties hold.
\begin{enumerate}
\item
The function $\bfe$ is regular with an infinite convergence radius,
and $\bfe(z,\theta)=E(z)$ for all $z\in\Omega$.

\item
The function $\bfe$ satisfies the following $\tau$-difference equation:
\begin{equation}\label{e:taudifference}
\bfe^{(2)}=\frac{1}{t-\theta^{q^2}}\bigl(g^q\bfe^{(1)}+\Delta\bfe\bigr).
\end{equation}

\item We have the following $u$-expansion, for $|t|<q^q$ and $|u(z)|$ sufficiently small:
\begin{equation*}\label{e:uexpansion}
\bfe(z)=u+\cdots\in u\FF_q[t,\theta][[u^{q-1}]].
\end{equation*}

\item
The function $\bfe$ satisfies the following functional equations,
for all $z\in\Omega$ and all
$\gamma\in \Gamma$:
\begin{equation}\label{ee:transformation}
\bfe(\gamma(z))=\frac{J_{\gamma}\bfj}{\det\gamma}\bigl(
\bfe(z)-\frac{1}{\widetilde{\pi}}\bfl\bigr).
\end{equation}
\end{enumerate}
\end{Proposition}

\noindent\emph{Proof.} Follows easily from Proposition~\ref{proph} (for the point 4 use (\ref{jmodularity})). \CVD

It follows at once from Propositions \ref{proph} and \ref{propE} that $\bfh$ and $\bfe$ are
almost-$A$-quasi-modular forms.

\begin{Corollaire}
we have
$$\bfh\in\widetilde{\mathcal{M}}^{\le 0}_{q,1,1} \qquad \text{and}\qquad\bfe\in\widetilde{\mathcal{M}}^{\le 1}_{1,1,1}.$$
\end{Corollaire}

\section{Proof of Theorem~\ref{estimate}}\label{theproof}

We prove here Theorem~\ref{estimate}, following the same method of \cite{Pe}.
We need a few more notations.
If $f:\Omega\rightarrow C$ is a non-zero Drinfeld quasi-modular form having a $u$-expansion
(for $|u(z)|$ small)
\begin{equation*}%\label{uexpansion}
f=\sum_{m\geq 0}f_m u^m,
\end{equation*}
we denote by $\nu_{\infty}(f)$ its {\em vanishing order at infinity}, that is
$$ \nu_{\infty}(f):=\min\{m\mid f_m\not=0\}. $$
We extend the notation $\nu_{\infty}(\bsb{f})$ to regular functions $\bsb{f}\in\bsb{\mathcal{O}}$ as follows.
If $\bsb{f}\in\bsb{\mathcal{O}}$, non-zero, then by definition %(see Definition~\ref{defiregular})
there exists $c>0$ such that, for $(t_0,z)\in C\times\Omega$ with $|t_0|<c$ and $|u(z)|<c$, one has
\begin{equation*}
f(z,t_0)=\sum_{m\geq 0}f_m(t_0) u^m,
\end{equation*}
where $f_m=\sum_{n\geq 0}f_{nm} t^n\in C[[t]]$ is a formal series. We define
$$ \nu_{\infty}(\bsb{f}):=\min\{m\mid f_m\not=0\}. $$

Following \cite{Pe}, we introduce the $\TT_{>0}$-algebra
$$ \MM^\dag=\TT_{>0}[g,h,\bsb{E},\bsb{E}^{(1)}]=\TT_{>0}[g,h,\bsb{E},\bsb{h}], $$
the equality being an easy consequence of the definition of $\bfe$ and of (\ref{taudifferenceh})
(see also \cite[Lemma~16]{Pe}). The functions $g,h,\bsb{E},\bsb{h}$ are algebraically independent over
$\TT_{>0}$ by \cite[Proposition~14]{Pe}.
For $\mu,\nu\in\ZZ$, $m\in\ZZ/(q-1)\ZZ$, we denote by
$\MM^{\dag}_{\mu,\nu,m}$ the sub-$\TT_{>0}$-module of $\MM^{\dag}$ consisting of
forms of weight $(\mu,\nu)$ and type $m$.

By (\ref{e:taudifference}), it is clear that $\MM^\dag$ is stable under twisting. More precisely,
we have the following result, already remarked in \cite{Pe}:

\begin{Lemma}\label{zzzz}
Let $\mu,\nu\in\ZZ$, $m\in\ZZ/(q-1)\ZZ$ and $k\in\ZZ_{\ge 0}$.
If $\bsb{f}\in \MM^{\dag}_{\mu,\nu,m}$, then $\bsb{f}^{(k)}\in \MM^{\dag}_{q^k\mu,\nu,m}$
\end{Lemma}

\noindent {\em Proof}. It suffices to prove the assertion for $k=1$ and for $\bsb{f}\in\{g,h,\bfh,\bfe\}$.
If $\bsb{f}=g$ or $\bsb{f}=h$ the result is clear. If $\bsb{f}=\bfh$ this is clear too since
$\bfh^{(1)}=\Delta\bfe\in \MM^{\dag}_{q^2,1,1}$. Suppose now that $\bsb{f}=\bfe$.
We have $\bfe^{(1)}=\frac{\bfh^{(2)}}{\Delta^q}$ by definition of $\bfe$. But then
the $\tau$-difference equation (\ref{taudifferenceh}) implies immediately $\bfe^{(1)}\in \MM^{\dag}_{q,1,1}$. \CVD

We recall the following multiplicity estimate for elements of $\MM^{\dag}_{\mu,\nu,m}$:

\begin{Lemma}\label{estimatedag}
Let $\bsb{f}$ be a non-zero element of $\MM^{\dag}_{\mu,\nu,m}$ with $\nu\not=0$. Then
\begin{equation*}
\nu_{\infty}(\bsb{f})\le\mu\nu.
\end{equation*} 
\end{Lemma}

\noindent {\em Proof.} \cite[Proposition~19]{Pe}. \CVD

The following lemma is a refinement of \cite[Lemma~24]{Pe}.

\begin{Lemma}\label{lemma1}
Let $\mu$, $\nu$ be two integers such that $\nu\geq 0$ and
$\mu-\nu\geq 6(q^2-1)$. Define $V:=\rg_{\TT_{>0}}(\MM^{\dag}_{\mu,\nu,m})$. Then
$$ \frac{(\mu-\nu)^2}{3(q^2-1)(q-1)}\leq V
\leq \frac{3}{(q^2-1)(q-1)}(\mu-\nu)^2.$$
\end{Lemma}

\noindent\textit{Proof.}
A basis of $\MM^{\dag}_{\mu,\nu,m}$ is given by
$$ (\phi_{is}\bfh^s\bfe^{\nu-s})_{0\le s\le\nu, 1\le i \le\sigma(s)},$$
where, for all $s$, $(\phi_{is})_{1\le i \le\sigma(s)}$ is a basis of $M_{\mu-s(q-1)-\nu,m-\nu}$.
Thus
$$ V=\sum_{0\leq s\leq\lfloor\frac{\mu-\nu}{q-1}\rfloor}
\dim_C(M_{\mu-s(q-1)-\nu,m-\nu}).$$
Now, we have (\eg\ \cite[Proposition 4.3]{cor}),
$$ \frac{\mu-\nu-s(q-1)}{q^2-1}-1\leq\dim_C(M_{\mu-s(q-1)-\nu,m-\nu})
\leq\frac{\mu-\nu-s(q-1)}{q^2-1}+1. $$
A simple computation gives
$$ (\lfloor x\rfloor+1)(\frac{1}{q+1}(x-\lfloor x\rfloor/2)-1) \leq V
\leq (\lfloor x\rfloor+1)(\frac{1}{q+1}(x-\lfloor x\rfloor/2)+1)$$
with $x=\frac{\mu-\nu}{q-1}$.
Using the inequalities $x/2\leq x-\lfloor x\rfloor/2\leq x/2+1/2$
and $x\leq \lfloor x\rfloor +1 \leq x+1$, we get
$$ x(\frac{x}{2(q+1)}-1)\leq V\leq (x+1)(\frac{x}{2(q+1)}+\frac{1}{2(q+1)}+1).$$
Since $ \frac{x}{2(q+1)}-1\geq \frac{x}{3(q+1)}$, $x+1\leq 2x$, $\frac{1}{2(q+1)}\leq \frac{x}{2(q+1)}$ and  $1\leq \frac{x}{2(q+1)}$,
we obtain the result. \CVD

\begin{Lemma}\label{lemma2}
Let $\mu$, $\nu$ be two integers such that
$\nu\geq 1$ and $\mu-\nu\geq 6(q^2-1)$. Let $m$ be an element of
$\{0,1,\ldots,q-2\}$, and suppose that $\MM^{\dag}_{\mu,\nu,m}\not=\{0\}$.
There exists a form $\bsb{f}_{\mu,\nu,m}\in \MM^{\dag}_{\mu,\nu,m}$, non-zero, satisfying
the following properties.

\begin{enumerate}
\item We have
\begin{equation}\label{minoration}
\nu_{\infty}(\bsb{f}_{\mu,\nu,m})\geq \frac{(\mu-\nu)^2}{9(q^2-1)}.
\end{equation}

\item The function $\bsb{f}_{\mu,\nu,m}$ has a $u$-expansion of the form
\begin{equation}\label{auxiliaire}
\bsb{f}_{\mu,\nu,m}=u^m\sum_{n\geq n_0}b_nu^{n(q-1)}
\end{equation}
with $b_{n_0}(t)\not=0$, $b_n(t)\in\FF_q[t,\theta]$ for all $n$, and
\begin{equation}\label{degbn0}
\deg_tb_{n_0}\leq \nu\log_q(\frac{3(\mu-\nu)^2}{2(q^2-1)(q-1)})+\nu\log_q(\frac{\mu\nu}{q-1}).
\end{equation}
\end{enumerate}
\end{Lemma}

\noindent {\em Proof.} Set $V:=\rg_{\TT_{>0}}(\MM^{\dag}_{\mu,\nu,m})$ as above, and
put $U=\lfloor V/2 \rfloor$. It follows from the proof of \cite[Proposition~5]{Pe}
that there exists a non-zero form $\bsb{f}_{\mu,\nu,m}\in \MM^{\dag}_{\mu,\nu,m}$
having an expansion of the form (\ref{auxiliaire}) with
$b_{n_0}(t)\not=0$, $b_n(t)\in\FF_q[t,\theta]$ for all $n$, and satisfying
\begin{equation}\label{eee}
U \le n_0\le \frac{\mu\nu}{q-1} \quad \text{and}\quad
\deg_tb_n\le \nu\log_qU+\nu\log_q\max\{1,n\}.
\end{equation}
We have obviously
\begin{equation}\label{mmm}
\nu_{\infty}(\bsb{f}_{\mu,\nu,m})=m+n_0(q-1)\ge n_0(q-1).
\end{equation}
We deduce from this
$$\nu_{\infty}(\bsb{f}_{\mu,\nu,m})\geq (q-1)U\geq \frac{(q-1)V}{3}, $$
hence the bound (\ref{minoration}) by Lemma~\ref{lemma1}.
The inequality (\ref{degbn0}) follows from (\ref{eee}), from the estimate
$U\le V/2$ and from Lemma~\ref{lemma1}. \CVD

\noindent {\em Proof of Theorem~\ref{estimate}.}

Let $f\in C[E,g,h]$ be a non-zero quasi-modular form of weight $w$, type $m$ and
depth $l\geq 1$.
Without loss of generality we may suppose that $f$ is irreducible.
\smallskip
We choose
$$\nu=1 $$
and
$$\mu=12(q^2-1)(w-l).$$
Since
$$ \mu-\nu\geq 12(q^2-1)-1\geq 6(q^2-1),$$
we may apply Lemma~\ref{lemma2} and we thus get the existence of a form
$\bsb{f}_{\mu,\nu,m}$.
Let $k$ be the smallest integer $\geq 0$ such that
\begin{equation}\label{choix_k}
q^k>3\log_q\mu.
\end{equation}
Using Lemma~\ref{lemma2} and its notation, we have
\begin{eqnarray*}
\deg_tb_{n_0}& \leq & \nu\log_q(\frac{3(\mu-\nu)^2}{2(q^2-1)(q-1)})+\nu\log_q(\frac{\mu\nu}{q-1})\cr
 & \leq & \log_q(\mu^2)+\log_q \mu=3\log_q\mu,
\end{eqnarray*}
hence $\deg_tb_{n_0}<q^k$.
We then define
$$f_k:=\left.\bsb{f}_{\mu,\nu,m}^{(k)}\right|_{t=\theta}.$$
By \cite[Lemma~22]{Pe}, we have $\nu_{\infty}(f_k)=\nu_{\infty}(\bsb{f}_{\mu,\nu,m}^{(k)})$,
hence
\begin{equation}\label{relationnu}
\nu_{\infty}(f_k)=q^k\nu_{\infty}(\bsb{f}_{\mu,\nu,m}).
\end{equation}
On the other hand, it follows from Lemma~\ref{zzzz} that
$$ w(f_k)=q^k\mu+\nu=q^k\mu+1.$$
Suppose first that $f$ does not divide $f_k$. Consider
the resultant $$\rho:=\res_E(f,f_k)\in C[g,h].$$ Since $\rho$ is a non-zero modular form,
we have $\nu_{\infty}(\rho)\le w(\rho)/(q+1)$. Thus, since there exist polynomials
$A,B\in C[E,g,h]$ with $Af+Bf_k=\rho$, we find
\begin{equation}\label{maj}
\min\{\nu_{\infty}(f);\nu_{\infty}(f_k)\}\leq\nu_{\infty}(\rho)
\leq \frac{w(\rho)}{q+1}=\frac{\nu(w-l)+q^k\mu l}{q+1}.
\end{equation}
If we had $\min\{\nu_{\infty}(f);\nu_{\infty}(f_k)\}=\nu_{\infty}(f_k)$,
then (\ref{maj}), (\ref{relationnu}) and (\ref{minoration}) would imply
$$ q^k\frac{(\mu-\nu)^2}{9(q^2-1)} \leq \frac{\nu(w-l)+q^k\mu l}{q+1},$$
or
\begin{equation}\label{cond}
q^k \frac{(\mu-1)^2}{9(q-1)}\leq w-l+q^k\mu l.
\end{equation}
But on one hand
\begin{equation*}
q^k \frac{(\mu-1)^2}{9(q-1)}> q^k\frac{\mu^2}{18(q-1)}\geq
\frac{3\mu\log_q\mu}{18(q-1)}\geq\frac{\mu}{6(q-1)}\geq 2(w-l),
\end{equation*}
and on the other hand
\begin{equation*}
q^k \frac{(\mu-1)^2}{9(q-1)}\geq q^k\frac{\mu^2}{18(q-1)}\geq
q^k\, \frac{12(q^2-1)\mu l}{18(q-1)} \geq 2 q^k\mu l.
\end{equation*}
This contradicts (\ref{cond}), hence
$\min\{\nu_{\infty}(f);\nu_{\infty}(f_k)\}=\nu_{\infty}(f)$, and
(\ref{maj}) gives the following upper bound:
\begin{equation}\label{majoration}
\nu_{\infty}(f) \leq \frac{(w-l)+q^k\mu l}{q+1}.
\end{equation}

Now, we have $\mu\geq 12(q^2-1)\geq q$, so $q^k>3$ by (\ref{choix_k})
and thus $k\geq 1$. It then follows from the definition of $k$ that
\begin{equation}\label{majqk}
q^{k-1}\leq 3 \log_q\mu.
\end{equation}
But the definition of $\mu$ gives
\begin{equation}\label{majmu}
\mu\leq 12 (q^2-1) (w-l),
\end{equation}
hence
\begin{equation}\label{majlogmu}
\log_q\mu\leq \log_q(12)+2+\log_q(w-l)\leq 7\max\{1,\log_q(w-l)\}.
\end{equation}
By (\ref{majqk}), (\ref{majlogmu}) and (\ref{majmu}) we get
$$ q^k\mu l \leq 252\, q(q^2-1) l (w-l)\max\{1,\log_q(w-l)\}.$$
The estimate $(\ref{majoration})$ now gives:
\begin{equation*}
\nu_{\infty}(f)\leq 252\, q^2 l(w-l)\max\{1,\log_q(w-l)\}
\end{equation*}
which implies the bound of the theorem.

Suppose now that $f$ divides $f_k$. Then
\begin{equation*}
\nu_{\infty}(f)\leq \nu_{\infty}(f_k)= q^k
\nu_{\infty}(\bsb{f}_{\mu,\nu,m})\leq q^k \mu \nu=q^k\mu
\end{equation*}
by Lemma~\ref{estimatedag}. Hence, again by (\ref{majqk}), (\ref{majlogmu}) and (\ref{majmu}):
$$ \nu_{\infty}(f)\leq 252\, q(q^2-1) (w-l)\max\{1,\log_q(w-l)\}.$$
\CVD

\section{$A$-modular and $A$-quasi-modular forms}\label{aformes}

Let us define, for $i=0,\ldots,q$:
$$\bsb{h}_i:=\widetilde{\pi}^i\bsb{s}_\carlitz^{-i}h\bsb{s}_2^i,$$ so that $\bsb{h}_0=h$ and $\bsb{h}_1=\bsb{h}$, the function introduced in \cite{Pe} and in
Section~\ref{preliminaries}. We recall that in
Section~\ref{preliminaries}, we have seen that the functions $E,g,h,\bsb{s}_2,\bsb{E},\bsb{h}$ are regular (and the radii of $E,g,h,\bsb{E},\bsb{h}$ are infinite).
It is a simple exercise to show that also the $\bsb{h}_i$'s
 (for $i=0,\ldots,q$) are regular, of infinite radius. 
 % If $\bsb{f}$ belongs to $\bsb{\mathcal{O}}$ and has radius $r$,
%then also $\bsb{f}^{(k)}$ belongs to $\bsb{\mathcal{O}}$ for all $k\geq 0$ and has radius at least [?] $r^{q^k}$.

We recall that $\bsb{s}_2$ is regular, with radius $q$.
We will use again that $\bsb{s}_2(z_0,t)\in C[[t]]^\times$ is a unit in the formal series, for all $z_0\in\Omega$ fixed.
At once, for all $t_0$ with $|t_0|$ small,
$\bsb{s}_2(z,t_0)$ identifies with a unit of $C[[u]]$ by (\ref{s2:uexpansion}).

\subsection{Almost $A$-quasi-modular forms}

Recall that we have defined the $\TT_{>0}$-algebra $\widetilde{\mathcal{M}}$ of
almost-$A$-quasi-modular forms in Section~\ref{aaaaa}. It contains the five algebraically
independent functions $E,g,h,\bsb{E},\bsb{E}^{(1)}$. However, in \cite{Pe}
there is no information about the structure of this algebra (the multiplicity estimate that was the main objective of that paper
only required the use of the fourth-dimensional algebra $\MM^\dag=\TT_{>0}[g,h,\bsb{E},\bsb{E}^{(1)}]$). In this paper we enlighten part of the structure of $\widetilde{\mathcal{M}}$.

We denote by $\mathcal{M}$ the $\TT_{>0}$-sub-algebra of $\widetilde{\mathcal{M}}$ generated by
%functions of $\widetilde{\mathcal{M}}$ which satisfy the functional equations
%(\ref{almostquasimodularforms}) with $\bsb{f}_{i,j}=0$ if $i+j>0$.
almost $A$-quasi-modular forms of depth $0$.
Obviously, $\mathcal{M}$ inherits the graduation from $\widetilde{\mathcal{M}}$ and we can write:
$$\mathcal{M}=\bigoplus_{\mu,\nu,m}\mathcal{M}_{\mu,\nu,m}$$ with
$\mathcal{M}_{\mu,\nu,m}=\widetilde{\mathcal{M}}_{\mu,\nu,m}^{\le 0}$.
%the appropriate $\TT_{>0}$-submodule of $\widetilde{\mathcal{M}}^{\leq l}_{\mu,\nu,m}$ (for all $l$). 
%Using the functional equations (\ref{twistss2}) it is easy
%to show that $\bsb{s}_2\in\mathcal{M}_{-1,1,0}$, for all $i=0,\ldots,q$, $\bsb{h}_i\in\mathcal{M}_{q+1-i,i,1}$,
%$\bsb{E}\in\widehat{\mathcal{M}}^{\leq 1}_{1,1,1}$ and $\bsb{E}^{(1)}\in\widehat{\mathcal{M}}^{\leq 1}_{q,1,1}$. It is also easy to show that
%$\bsb{E}^{(1)}$

We will prove the following theorem which supplies partial information on the structure of $\widetilde{\mathcal{M}}$:

\begin{Theorem}\label{theoremofstructure}
The $\ZZ^2\times\ZZ/(q-1)\ZZ$-algebra $\widetilde{M}$ has dimension five over $\TT_{>0}$ and we have 
\begin{equation}\label{equality1}
\widetilde{\mathcal{M}}=\mathcal{M}[E,\bsb{E}].
\end{equation}
Moreover,
the following inclusions hold, implying that $\mathcal{M}$ has dimension three over $\TT_{>0}$:
\begin{equation}\label{doubleinclusions}
\TT_{>0}[g,h,\bfs_2]\subset\mathcal{M}\subset \TT_{>0}[g,h,\bfs_2,\bfs_2^{-1}].
\end{equation}
\end{Theorem}

As exercises to familiarise with Theorem \ref{theoremofstructure}, the reader can verify the following properties (hint: do not forget to use, for example, the functional equations (\ref{twistss2})):
\begin{enumerate}
\item $M\subset\mathcal{M}$, $\widetilde{M}\subset\widetilde{\mathcal{M}}$.
\item $g\in\mathcal{M}_{q-1,0,0}$, $h\in\mathcal{M}_{q+1,0,1}$ and $E\in\widetilde{\mathcal{M}}_{2,0,1}^{\leq 1}$.
\item $\bsb{s}_2\in\mathcal{M}_{-1,1,0}$ but $\bsb{s}_2^{(1)}\not\in\widetilde{\mathcal{M}}$.
\item For all $i=0,\ldots,q$, $\bsb{h}_i\in\mathcal{M}_{q+1-i,i,1}$.
\item $\bsb{E}\in\widetilde{\mathcal{M}}^{\leq 1}_{1,1,1}$ and $\bsb{E}^{(1)}\in\widetilde{\mathcal{M}}^{\leq 1}_{q,1,1}$.
\item There exists an element $\lambda\in\FF_q(t,\theta)^\times$ (compute it!) such that
$\bsb{E}^{(1)}=\lambda(\bsb{h}_1+g\bsb{E})$.
\item $\MM^\dag=\TT_{>0}[g,h,\bsb{E},\bsb{E}^{(1)}]=\TT_{>0}[g,h,\bsb{E},\bsb{h}_1]$.
\item $\bsb{E}=-\sqrt[q]{\bsb{h}_q^{(1)}}$.
\end{enumerate}

\subsection{Structure of $\widetilde{\mathcal{M}}$}

In this subsection, we show that (\ref{equality1}) holds; this will be a consequence of Proposition \ref{widetildeMstructure} that will need the proposition
below.

\begin{Proposition}\label{leadingforms} For $l>0$,
let $\bsb{f}$ be an element of $\widetilde{\mathcal{M}}^{\leq l}_{\mu,\nu,m}\setminus\widetilde{\mathcal{M}}^{\leq l-1}_{\mu,\nu,m}$
satisfying 
(\ref{almostquasimodularforms}). Then, if $i+j=l$, we have, in (\ref{almostquasimodularforms}), $\bsb{f}_{i,j}\in\mathcal{M}_{\mu-2i-j,\nu-j,m-i-j}$.
\end{Proposition}
\noindent\emph{Proof.} Let us consider three matrices $\mathcal{A},\mathcal{B},\mathcal{C}\in\Gamma$ as follows:
\begin{equation}\label{ABC}
\mathcal{A}=\begin{pmatrix}a & b \\ c & d\end{pmatrix},\quad \mathcal{B}=\begin{pmatrix}\alpha & \beta \\ \gamma & \delta\end{pmatrix},\quad \mathcal{C}={\mathcal{A}}\cdot{\mathcal{B}}
=\begin{pmatrix}* & * \\ x & y\end{pmatrix}.
\end{equation}
We recall from \cite{Pe} that 
we have the functional equations:
\begin{eqnarray}
L_{\mathcal{A}}(\mathcal{B}(z))&=&\det(\mathcal{B})^{-1}J_{\mathcal{B}}(z)^2(L_{\mathcal{C}}(z)-L_{\mathcal{B}}(z)),\nonumber\\
\boldsymbol{L}_{\mathcal{A}}(\mathcal{B}(z))&=&\det(\mathcal{B})^{-1}J_{\mathcal{B}}(z)\boldsymbol{J}_{\mathcal{B}}(z)(\boldsymbol{L}_{\mathcal{C}}(z)-\boldsymbol{L}_{\mathcal{B}}(z)).\label{eqLL}
\end{eqnarray}

We now compute, by using (\ref{eqLL}) and the fact that $(\bsb{J}_{\gamma})_{\gamma\in\Gamma}$ is a factor of automorphy (cf. (\ref{factorofautomorphy})):
\begin{eqnarray*}
\lefteqn{ \bsb{f}(\mathcal{C}(z))=\bsb{f}(\mathcal{A}(\mathcal{B}(z)))=}\\
&=&\det(\mathcal{A})^{-m}J_{\mathcal{A}}(\mathcal{B}(z))^\mu\boldsymbol{J}_{\mathcal{A}}(\mathcal{B}(z))^\nu
\sum_{i+j\leq l}\bsb{f}_{i,j}(\mathcal{B}(z))L_{\mathcal{A}}(\mathcal{B}(z))^i\boldsymbol{L}_{\mathcal{A}}(\mathcal{B}(z))^j\\
&=&\det(\mathcal{A})^{-m}J_{\mathcal{A}}(\mathcal{B}(z))^\mu\boldsymbol{J}_{\mathcal{A}}(\mathcal{B}(z))^\nu\times\\
& &\sum_{i+j\leq l}\bsb{f}_{i,j}(\mathcal{B}(z))(\det(\mathcal{B})^{-1}J_{\mathcal{B}}(z)^2(L_{\mathcal{C}}(z)-L_{\mathcal{B}}(z)))^i(\det(\mathcal{B})^{-1}J_{\mathcal{B}}(z)\times\\
& &\boldsymbol{J}_{\mathcal{B}}(z)(\boldsymbol{L}_{\mathcal{C}}(z)-\boldsymbol{L}_{\mathcal{B}}(z)))^j\\
&=&\det(\mathcal{A})^{-m}J_{\mathcal{A}}(\mathcal{B}(z))^\mu\boldsymbol{J}_{\mathcal{A}}(\mathcal{B}(z))^\nu\times\\
& &\sum_{i+j\leq l}\bsb{f}_{i,j}(\mathcal{B}(z))\det(\mathcal{B})^{-i-j}J_{\mathcal{B}}(z)^{2i+j}\boldsymbol{J}_{\mathcal{B}}(z)^j\times\\
& &\sum_{s=0}^{i}(-1)^{i-s}\binom{i}{s}L_{\mathcal{C}}^sL_{\mathcal{B}}^{i-s}
\sum_{s'=0}^j(-1)^{j-s'}\binom{j}{s'}\boldsymbol{L}_{\mathcal{C}}^{s'}\boldsymbol{L}_{\mathcal{B}}^{j-s'}\\
&=&\det(\mathcal{A})^{-m}J_{\mathcal{A}}(\mathcal{B}(z))^\mu\boldsymbol{J}_{\mathcal{A}}(\mathcal{B}(z))^\nu\times\\
& &\sum_{s+s'\leq l}L_{\mathcal{C}}^{s}\boldsymbol{L}_{\mathcal{C}}^{s'}\sum_{i=s}^l\times\\ & &\sum_{j=s'}^{l-i}\bsb{f}_{i,j}(\mathcal{B}(z))\det(\mathcal{B})^{-i-j}
J_{\mathcal{B}}^{2i+j}\boldsymbol{J}_{\mathcal{B}}^{j}(-1)^{i+j-s-s'}\binom{i}{s}\binom{j}{s'}L_{\mathcal{B}}^{i-s}\boldsymbol{L}_{\mathcal{B}}^{j-s}.
\end{eqnarray*}
On the other side, we see that
\begin{eqnarray*}
\bsb{f}(\mathcal{C}(z))&=&\det(\mathcal{A}\mathcal{B})^{-m}J_\mathcal{C}^\mu\boldsymbol{J}_\mathcal{C}^\nu\sum_{i+j\leq l}\bsb{f}_{i,j}L_\mathcal{C}^i\boldsymbol{L}_{\mathcal{C}}^j\\
&=&\det(\mathcal{A}\mathcal{B})^{-m}J_\mathcal{A}(\mathcal{B}(z))^\mu J_\mathcal{B}^\mu\boldsymbol{J}_\mathcal{A}(\mathcal{B}(z))^\nu\boldsymbol{J}_\mathcal{B}^\nu\sum_{i+j\leq l}\bsb{f}_{i,j}L_\mathcal{C}^i\boldsymbol{L}_{\mathcal{C}}^j.
\end{eqnarray*}
We let $\mathcal{A}$ span $\mathbf{GL}_2(A)$ leaving $\mathcal{B}$ fixed at the same time.
Since in this way the matrix $\mathcal{C}$ covers the whole group $\Gamma$,
if $s,s'$ are such that $s+s'\leq l$, then
\begin{eqnarray*}
\lefteqn{\det(\mathcal{B})^{-m}J_\mathcal{B}^\mu\boldsymbol{J}_\mathcal{B}^\nu \bsb{f}_{s,s'}(z)=}\\ &=&\sum_{i=s}^l\sum_{j=s'}^{l-i}\bsb{f}_{i,j}(\mathcal{B}(z))\det(\mathcal{B})^{-i-j}J_\mathcal{B}^{2i+j}\boldsymbol{J}_{\mathcal{B}}^{j}(-1)^{i+j-s-s'}\binom{i}{s}\binom{j}{s'}L_{\mathcal{B}}^{i-s}\boldsymbol{L}_{\mathcal{B}}^{j-s'}.
\end{eqnarray*}
If $s+s'=l$, there is only one term in the sum above, corresponding to $i=s,j=s'$. Thus we find
%and the sum above is non-zero, we have in it the only non-zero contribution for $i=s,j=s'$. Therefore $\bsb{f}_{s,s'}\neq 0$ and
$$\bsb{f}_{s,s'}(\mathcal{B}(z))=\det(\mathcal{B})^{s+s'-m}J_\mathcal{B}^{\mu-2s-s'}\boldsymbol{J}_\mathcal{B}^{\nu-s'}\bsb{f}_{s,s'}(z).$$
This identity of formal series of $C[[t]]$ holds for every $\mathcal{B}\in\Gamma$ and every $z\in\Omega$.
Since on the other side we know already that $\bsb{f}_{s,s'}\in\boldsymbol{\mathcal{O}}$,
the proposition follows.\CVD

\begin{Proposition}\label{widetildeMstructure}
Let $\bsb{f}$ be in $\widetilde{\mathcal{M}}^{\leq l}_{\mu,\nu,m}$.
Then $$\bsb{f}\in\bigoplus_{i+j\leq l}E^i\boldsymbol{E}^j\mathcal{M}_{\mu-2i-j,\nu-j,m-i-j}.$$
\end{Proposition}
\noindent\emph{Proof.}
Let us assume that $\bsb{f}$ is an almost $A$-quasi-modular form, such that for $\gamma\in\Gamma$,
(\ref{almostquasimodularforms}) holds.
By Proposition \ref{leadingforms}, if $i+j=l$ in (\ref{almostquasimodularforms}), then
$\bsb{f}_{i,j}\in\mathcal{M}_{\mu-2i-j,\nu-j,m-i-j}$.
By the functional equations of $E,\bsb{E}$ (\ref{e:transformation}) and (\ref{ee:transformation}),
the almost $A$-quasi-modular form $\bsb{\rho}_{i,j}:=E^i\boldsymbol{E}^j$, for $\gamma\in\Gamma$, transforms like:
$$\bsb{\rho}_{i,j}(\gamma(z))=\det(\gamma)^{-i-j}J_\gamma^{2i+j}\boldsymbol{J}_\gamma^{j}\left(\bsb{\rho}_{i,j}(z)+\cdots+(-1)^{i+j}\frac{1}{\widetilde{\pi}^{i+j}}L_\gamma^i\boldsymbol{L}_\gamma^j\right).$$
Hence, the function:
$$\bsb{f}':=\bsb{f}-(-1)^l\widetilde{\pi}^l\sum_{i+j=l}\bsb{\rho}_{i,j}\bsb{f}_{i,j}$$
is an almost $A$-quasi-modular form with same weight and type as $\bsb{f}$, and depth strictly less than $l$.
We can apply Proposition \ref{leadingforms} again on this form, and construct in this way 
another almost $A$-quasi-modular form of depth strictly less than $l-1$ and so on. Since the depth is positive,
we will end the inductive process with an almost $A$-quasi-modular form of depth $\leq 0$.
Summing up all the terms, we obtain what we wanted, noticing that the proposition immediately implies that (\ref{equality1}) holds.
\CVD

\subsection{The algebra $\mathcal{M}$}\label{pasdidee}

We prove (\ref{doubleinclusions}) in this subsection. We first need some preliminaries.

\begin{Lemma}\label{lemma:nozeros}
Let $(z_0,t_0)$ be in $\Omega\times (C\setminus \{\theta,\theta^q,\theta^{q^2},\ldots\})$.
There exists $\gamma\in\Gamma$ such that $\bfs_2(\gamma(z_0),t_0)\not=0$.
\end{Lemma}
\noindent\emph{Proof.} 
We can suppose that $(z_0,t_0)$ are such that $\bfs_2(z_0,t_0)=0$, otherwise the result is trivial. 

If $\gamma=\begin{pmatrix}a & b \\ c & d\end{pmatrix}\in\Gamma$,  by
(\ref{eqvectorial}),
$$\bfs_2(\gamma(z_0),t_0)=(cz_0+d)^{-1}(\overline{c}(t_0)\bfs_1(z_0,t_0)+\overline{d}(t_0)\bfs_2(z_0,t_0)).$$

Let us choose $\gamma$ so that $c\in A$ satisfies $\overline{c}(t_0)\not=0$. If also $\bfs_2(\gamma(z_0),t_0)$ vanishes, then, $\bfs_1(z_0,t_0)=0$.
By (\ref{dethPsi}), $\bfs_\carlitz(t_0)=0$ (the form $h$ is holomorphic with no zeros on $\Omega$ so its inverse is holomorphic and does not vanish).
But the formula (\ref{scarlitz-product}) is contradictory with the vanishing of $\bfs_\carlitz(t_0)$ just obtained.\CVD

\begin{Remarque}{\em 
Thanks to Lemma \ref{lemma:nozeros}, we will use $\bsb{s}_2$ to
detect the structure of our automorphic functions following closely the usual procedure for Drinfeld modular forms, where one uses the fact that
$h$ does not vanish in $\Omega$ to deduce the structure of $S_{w,m}$ (space of cusp forms) from the structure of $M_{w-q-1,m-1}$.
}\end{Remarque}

We will also need the following proposition:
\begin{Proposition}\label{proposition:mathcalM} For given $\mu,\nu,m$,
let $\bsb{f}$ be a non-zero element of $\mathcal{M}_{\mu,\nu,m}$.
Then, we have $\nu\geq -\mu$ and the function:
$$F:=\bsb{f}\bfs_2^{-\nu}$$ belongs to $M_{\mu+\nu,m}\otimes_{C} \TT_{>0}$. 
\end{Proposition}
\noindent\emph{Proof.} By (\ref{s2modularity}) the function $F$ satisfies, for all $z\in\Omega$, the following identities
in $C[[t]]$:
$$F(\gamma(z))=\det(\gamma)^{-m}J_\gamma^{\mu+\nu}F(z).$$
Let $t_0$ be such that $|t_0|<r$, with $r$ the minimum of the radius of $\bsb{f}$ and $q$, the radius of $\bsb{s}_2$. We know that
the functions $z\mapsto\bsb{f}(z,t_0)$ and $z\mapsto\bfs_2(z,t_0)$ are holomorphic functions $\Omega\rightarrow C$.
Therefore, $\varphi_{t_0}(z):=F(z,t_0)$ is a meromorphic function over $\Omega$ (more precisely, it is holomorphic if $\nu\leq 0$). 
Moreover, it is plain that for all $\gamma\in\Gamma$ and $z\in\Omega$ at which
$\varphi_{t_0}$ is defined, 
$$\varphi_{t_0}(\gamma(z))=\det(\gamma)^{-m}J_\gamma^{\mu+\nu}\varphi_{t_0}(z).$$ We also 
see that $\varphi_{t_0}$ is holomorphic at infinity. Indeed, for all $t_0$ with $|t_0|$ small, $\bsb{s}_2(z,t_0)$ is a unit of $C[[u]]$, by (\ref{s2:uexpansion}).
Moreover, there exists $c>0$ depending on $t_0$ such that $\bsb{s}_2(z,t_0)$ has no zeroes for $z$ such that $|u|=|u(z)|<c$ and this property 
is obviously shared with $\bsb{s}_2(z,t_0)^{-1}$.
Hence, $\varphi_{t_0}$ is a ``meromorphic Drinfeld modular form" of weight $\mu+\nu$ and type $m$ which is holomorphic at infinity.

We claim that in fact, $\varphi_{t_0}$ is holomorphic on $\Omega$ (hence regular) regardless of the sign of $\nu$ (this implies that it is a ``genuine" Drinfeld modular form).
Indeed, if the claim is false, $\nu>0$ and $\varphi_{t_0}$ has a pole at $z_0\in\Omega$,
as well as at any point $\gamma(z_0)$ with $\gamma\in\Gamma$.
Since $z\mapsto \bsb{f}(z,t_0)$ is holomorphic in $\Omega$, $\gamma(z_0)$ is a zero of $z\mapsto\bfs_2(z,t_0)$ for
all $\gamma\in\Gamma$ as above. But this is contradictory with Lemma \ref{lemma:nozeros}.

%Hence, we have at once $F\in M_{\mu+\nu,m}\otimes_C(C-\text{algebra of maps }B(0,r)\rightarrow C)$ (where $B(0,r)$ is the ``open" ball of radius $r$
%of center $0$) and $F\in \bsb{\mathcal{O}}$.

Hence, for all $t_0$ as above, $\varphi_{t_0}=F(\cdot,t_0)\in M_{\mu+\nu,m}$.
Let $(b_1,\ldots,b_\ell)$ be a basis of $M_{\mu+\nu,m}$ over $C$. 
It is easy to show, looking at the $(u_1,\ldots,u_\ell)$-expansion in $C[[u_1,\ldots,u_\ell]]$ of the function $\det(b_i(z_j))_{1\leq i,j\leq \ell}$
and the existence of the ``triangular" basis $(g^ih^j)$ in $M_{\mu+\nu,m}$, that there exist (obviously distinct)
$z_1,\ldots,z_\ell\in\Omega$ such that the matrix $(b_i(z_j))_{1\leq i,j\leq \ell}\in\mathbf{Mat}_{\ell\times\ell}(C)$ is non-singular.

For all $t$ with $|t|<r$, there exist
$d_1(t),\ldots,d_\ell(t)\in C$ such that 
$$F(z,t)=d_1(t)b_1(z)+\cdots+d_\ell(t)b_\ell(z).$$

Since for all $i=1,\ldots,\ell$, $F(z_i,t)$  belongs to $ \TT_{>0}$, we find:
$$d_1(t)b_1(z_j)+\cdots+d_\ell(t)b_\ell(z_j)\in  \TT_{>0},\quad j=1,\ldots,\ell.$$
Solving a linear system in the indeterminates $d_1,\ldots,d_\ell$ we find $d_1,\ldots,d_\ell\in  \TT_{>0}$.
Hence, $F\in M_{\mu+\nu,m}\otimes_{C}  \TT_{>0}$ and $\bsb{f}$ being non-zero, we have $\mu+\nu\geq 0$, concluding the proof of the proposition.\CVD

We deduce from  Proposition \ref{proposition:mathcalM}  and from the discussion above, the inclusions (\ref{doubleinclusions}).

\begin{Remarque}
{\em If $f$ is a Drinfeld modular form in $M_{w,m}$, $\bsb{f}_i:=\widetilde{\pi}^i\bsb{s}_{\carlitz}^{-i}f\bsb{s}_2^i$ belongs to
$\mathcal{M}_{w-i,i,m}$ for all $i\geq 0$, has infinite radius, and is such that $\bsb{f}_i|_{t=\theta}=f$ (this is how we have defined the forms $\bsb{h}_i$). 
For example, if $f=1$, we have $\bsb{1}_i=\widetilde{\pi}^i\bsb{s}_{\carlitz}^{-i}\bsb{s}_2^i$. 
Proposition \ref{proposition:mathcalM} gives a converse of this fact: every element of $\mathcal{M}_{\mu,\nu,m}$ comes from a Drinfeld 
modular form in this way. However, 
it is an open question to prove or disprove that $\bsb{s}_2^{-1}$ is regular. In fact,
{\em we expect that} $\bsb{s}_2^{-1}\not\in\bsb{\mathcal{O}}$.
Several arguments led us to believe that there exists a sequence $(z_n,t_n)_{n\geq 0}$ in $\Omega\times C$ with 
$\lim_{n\rightarrow\infty}t_n=0$ and $t_n\neq 0$ for all $n$,
such that $\bsb{s}_2(z_n,t_n)=0$ for all $n$ but, at the time being, 
we do not know how to precisely locate the zeroes of the function $(z,t)\mapsto\bsb{s}_2(z,t)$. }
\end{Remarque}

\begin{Question}\label{quest}Find the exact image of $\mathcal{M}$ in $\TT_{>0}[g,h,\bfs_2,\bfs_2^{-1}]$.\end{Question}

The following conjecture agrees with our guess that $\bsb{s}_2^{-1}$ is not regular.

\begin{Conjecture}\label{conj}We have $\mathcal{M}=\TT_{>0}[g,h,\bfs_2]$.\end{Conjecture}

In the next subsection we will construct a three-dimensional sub-algebra $\MM$ of  $\mathcal{M}$ of {\em $A$-modular forms} with a very precise structure, which also is finitely generated over $\TT_{>0}$.
To define it, we will use the Frobenius structure over $\bsb{\mathcal{O}}$. With this task in view, we will first
need to introduce the modules $\widehat{\mathcal{M}}^{\leq l}_{\mu,\nu,m}$ at the beginning of the next subsection.

\subsection{$A$-modular forms and their structure}

For $(\mu,\nu,m)$ an element of $\ZZ^2\times \ZZ/(q-1)\ZZ$ and $l\in\NN$,
we will use, in all the following, the sub-$\TT_{>0}$-module $\widehat{\mathcal{M}}^{\leq l}_{\mu,\nu,m}\subset\widetilde{\mathcal{M}}^{\leq l}_{\mu,\nu,m}$ whose elements $\bsb{f}$ satisfy functional equations of the more particular form:
$$\bsb{f}(\gamma(z),t)=\det(\gamma)^{-m}J_\gamma^\mu\bsb{J}_\gamma^\nu\left(\sum_{i\leq l}\bsb{f}_{i}\bsb{L}_\gamma^i\right),$$
with $\bsb{f}_{i}\in\bsb{\mathcal{O}}$ for all $i$. In such a kind of functional equation, there is no dependence on $L_\gamma$. Hence,  $E\not\in\widehat{\mathcal{M}}$ where 
$\widehat{\mathcal{M}}$ is the $\TT_{>0}$-algebra, graded by $\ZZ^2\times \ZZ/(q-1)\ZZ$, generated by the modules 
$\widehat{\mathcal{M}}^{\leq l}_{\mu,\nu,m}$. We deduce from the last subsection that $\widehat{\mathcal{M}}=\widetilde{\mathcal{M}}[\bsb{E}]$ is a $\TT_{>0}$-algebra of 
dimension $4$.

\begin{Definition}
{\em Let $\bsb{f}$ be an element of $\widehat{\mathcal{M}}^{\leq l}_{\mu,\nu,m}$.
We say that $\bsb{f}$ is an {\em $A$-quasi-modular form} of weight $(\mu,\nu)$, type $m$ and depth $\leq l$, if, for all $k\geq 0$, 
\begin{equation}\label{conditionforquasimodular}
\bsb{f}^{(k)}\in\widehat{\mathcal{M}}.
\end{equation}
If $l=0$, we will say that $\bsb{f}$ is an {\em $A$-modular form} of weight $(\mu,\nu)$ and type $m$.}
\end{Definition}

For $\mu,\nu\in\ZZ,m\in\ZZ/(q-1)\ZZ$, $l\in\ZZ_{\geq 0}$, we denote by $\widetilde{\MM}_{\mu,\nu,m}^{\leq l}$ the $\TT_{>0}$-module of
$A$-quasi-modular forms of weight $(\mu,\nu)$, type $m$ and depth $\leq l$. Again,
we have $$\widetilde{\MM}^{\leq l}_{\mu,\nu,m}\widetilde{\MM}^{\leq l'}_{\mu',\nu',m'}\subset\widetilde{\MM}^{\leq l+l'}_{\mu+\mu',\nu+\nu',m+m'}$$
and the $\TT_{>0}$-algebra $\widetilde{\MM}$ generated by all these modules is again graded by $\ZZ^2\times \ZZ/(q-1)\ZZ$ (and filtered by the depths).
The graded sub-algebra of $\widetilde{\MM}$ generated by $A$-modular forms is denoted by $\MM$. We have thus $\MM=\oplus_{\mu,\nu,m}\MM_{\mu,\nu,m}$,
where $\MM_{\mu,\nu,m}=\widetilde{\MM}^{\leq 0}_{\mu,\nu,m}$.

Earlier in this paper and in \cite{Pe}, we have introduced the $\TT_{>0}$-algebra $\MM^\dag=\TT_{>0}[g,h,\bsb{E},\bsb{E}^{(1)}]$, which turned out to be equal to
$\TT_{>0}[g,h,\bsb{E},\bsb{h}_1]$. We have $\MM^\dag\subset\widetilde{\MM}$ but these algebras are not equal. Indeed,
with the help of (\ref{twistss2}) (\footnote{Identity (\ref{twistss2}) also certifies that $\bsb{s}_2$ does not belong to $\widetilde{\MM}$.
Indeed, $\bsb{s}_2^{(1)}\not\in\widetilde{\mathcal{M}}$.}), one sees that $\TT_{>0}[g,h,\bsb{h}_1,\ldots,\bsb{h}_q]\subset\MM$.
The next theorem provides the equality of the latter two algebras. We notice that at the time being, we do not have a complete structure
theorem for the algebra $\widetilde{\MM}$.

\begin{Theorem}\label{structureofaforms} The algebra 
$\MM$ of $A$-modular forms is finitely generated of dimension three. Moreover, 
we have $\MM=\TT_{>0}[g,\bsb{h}_0,\bsb{h}_1,\ldots,\bsb{h}_q]$.
\end{Theorem}

\noindent\emph{Proof.} 
Since it is clear that 
$\TT_{>0}[g,\bsb{h}_0,\bsb{h}_1,\ldots,\bsb{h}_q]$ has dimension three, the first 
part of the statement is a consequence of the second that we prove now.
Let $\bsb{f}$ be a non-constant element of $\MM_{\mu,\nu,m}$. Since $\MM_{\mu,\nu,m}\subset\mathcal{M}_{\mu,\nu,m}$,
Proposition
\ref{proposition:mathcalM} implies that $\bsb{f}\bsb{s}_2^{-\nu}\in\TT_{>0}[g,h]$. Hence, there exists a non-zero 
element $\varphi\in M_{\mu+\nu,m}\otimes_{C}\TT_{>0}$ such that $\bsb{f}=\varphi\bsb{s}_2^\nu$.

Let $s$ be the integer $\nu_\infty(\bsb{f})$; we have that $s=\nu_\infty(\varphi)$ because $\nu_\infty(\bsb{s}_2)=0$ as we remarked earlier. We claim that 
$0\leq\nu\leq sq$. We first proceed to prove that $\nu\geq 0$.

Assume conversely that $\nu<0$.
The hypotheses and (\ref{twistss2}) imply, for $k\geq 0$, for all $z\in\Omega$ and for all $\gamma\in\Gamma$:
\begin{eqnarray*}
(\varphi\bsb{s}_2^{\nu})^{(k)}(\gamma(z),t)&=&J_\gamma^{\mu q^k}\bsb{J}_\gamma^{\nu}\det(\gamma)^{-m}\varphi(z)^{(k)}(\bsb{s}_2^{(k)}+
\alpha_k\bsb{L}_\gamma)^{\nu}\\
&=&\sum_{a,b,c}J_\gamma^a\bsb{J}_\gamma^b\det(\gamma)^{-c}\sum_j\bsb{f}_{a,b,c,j,k}\bsb{L}_\gamma^j,
\end{eqnarray*}
where $$\alpha_k=(t-\theta)
\frac{g^{\star}_{k-1}(z,t)\bfs_{\carlitz}(t)}{\widetilde{\pi}^{q^k+1}h(z)^{q^{k-1}}}$$ is the (non-zero) coefficient of $\bsb{L}_\gamma$ in the inner bracket of the right-hand side of (\ref{twistss2}), where all the sums have finitely many terms,
and where the $\bsb{f}_{a,b,c,j,k}$'s are regular functions. Comparing the two right-hand sides we get:
\begin{eqnarray*}
1&=&\det(\gamma)^{m}J_\gamma^{-\mu q^k}\bsb{J}_\gamma^{-\nu}(\varphi^{-1})^{(k)}\times\\
& &(\bsb{s}_2^{(k)}+\alpha_k\bsb{L}_\gamma)^{-\nu}\sum_{a,b,c}J_\gamma^a\bsb{J}_\gamma^b\det(\gamma)^{-c}\sum_j\bsb{f}_{a,b,c,j,k}\bsb{L}_\gamma^j
\end{eqnarray*}
which provides a contradiction with the non-vanishing of $\alpha_k$ because the constant function $z\mapsto 1$ is 
$A$-quasi-modular of weight $(0,0)$, type $0$ and depth $\leq 0$.

We have proved that $\nu\geq 0$. We now prove the upper bound $\nu\leq qs$. 
We can write $\varphi=h^s\psi$ with $\psi\in\TT_{>0}[g,h]$ such that $\nu_\infty(\psi)=0$. Since $\bsb{f}=h^s\psi\bsb{s}_2^\nu$,
we have, for all $k$, $\bsb{f}^{(k)}=h^{sq^k}\psi^{(k)}(\bsb{s}_2^\nu)^{(k)}$ and
$$\bsb{f}^{(k)}(\gamma(z))=\det(\gamma)^{-m}J_\gamma^{\mu q^k}\bsb{J}_\gamma^{\nu}h(z)^{sq^k}\psi^{(k)}(z)(\bsb{s}_2^{(k)}(z)+\alpha_k\bsb{L}_\gamma)^{\nu}.$$

A necessary condition for $\bsb{f}^{(k)}$ to be in $\widehat{\mathcal{M}}$ is that $h^{sq^k}\alpha_k^\nu$ belongs to $\bsb{\mathcal{O}}$ for all $k\geq 1$. 
Since $\nu_\infty(\alpha_k)=-q^{k-1}$ (when $k\geq 1$), this condition is equivalent to $sq-\nu\geq 0$. This gives
$\nu\leq qs$ as claimed.

Write now $\nu=qa+r$ with $0\leq r<q$ (Euclidean division), from which we deduce
$a\leq s$ and $a<s$ if $r\neq0$.

If $r=0$ then we have $\nu=qa$ and $h^s\bsb{s}_2^\nu=h^s\bsb{s}_2^{qa}=h^{s-a}h^a\bsb{s}_2^{qa}$ which belongs to
$\TT_{>0}\bsb{h}_0^{s-a}\bsb{h}_q^a$. If $r\not= 0$, then $s-a-1\geq 0$
and
$$ h^s\bsb{s}_2^\nu=h^{s-a-1}.h\bsb{s}_2^r.h^a\bsb{s}_2^{qa}\in \TT_{>0}\bsb{h}_0^{s-a-1}\bsb{h}_r\bsb{h}_q^a. $$
In all cases, we thus have $h^s\bsb{s}_2^\nu\in\TT_{>0}[\bsb{h}_0,\bsb{h}_1,\ldots,\bsb{h}_q]$.
Since $\bsb{f}=h^s\psi\bsb{s}_2^\nu$ with $\psi\in \TT_{>0}[g,h]$, we finally obtain
$\bsb{f}\in\TT_{>0}[g,\bsb{h}_0,\bsb{h}_1,\ldots,\bsb{h}_q]$ as required.\CVD

\section{Final remarks}

\begin{Remarque}\label{firstremark}{\em 
The image of the map (Anderson's twist) $\tau:\widetilde{\mathcal{M}}\rightarrow\bsb{\mathcal{O}}$ is not contained in $\widetilde{\mathcal{M}}$. At least, it can be proved
that given
$\bsb{f}\in \widetilde{\mathcal{M}}_{\mu,\nu,m}^{\leq l}$, we have $\tau^{k}\bsb{f}\in \widetilde{\mathcal{M}}[1/h]_{q^k\mu,\nu,m}^{\leq l}$ for all $k\geq 0$.

In Lemma \ref{zzzz} (cf. \cite{Pe}), we have showed that $\FF_q[[t]]$-linear Anderson's operator $\tau:\bsb{\mathcal{O}}\rightarrow\bsb{\mathcal{O}}$, which operates
on double formal series of $C[[t,u]]$ as:
$$\tau:\sum_{i,j\geq 0}c_{i,j}t^iu^j\mapsto \sum_{i,j\geq 0}c_{i,j}^qt^iu^{qj}$$
(the $c_{i,j}$'s being elements of $C$) defines an operator $\MM^\dag\rightarrow\MM^\dag$ and is $\TT_{>0}$-linear
on the modules $\MM^\dag_{\mu,\nu,m}\rightarrow\MM^\dag_{q\mu,\nu,m}$. Thanks to this result, it is possible to prove that every 
form of $\MM^\dag_{\mu,\nu,m}$ satisfies a non-trivial linear $\tau$-difference equation.

Since we still do not know the exact structure of the algebra $\widetilde{\MM}$, we presently cannot extend this observation to the settings of this paper.
}
\end{Remarque}

%\begin{Remarque}{\em 
%The $\FF_q$-linear Frobenius $F$ on $C[[t,u]]$ (which sends $x$ to $x^q$) splits as a product
%\begin{equation}\label{splitting}
%F=\tau\chi=\chi\tau,\end{equation}
%where $\chi$ is the $C[[u]]$-linear operator {\em adjoint} to the operator $\tau$, defined by:
%$$\chi:\sum_{i,j\geq 0}c_{i,j}t^iu^j\mapsto \tau\sum_{i,j\geq 0}c_{i,j}t^{qi}u^{j}.$$ It turns out, in a way which is complementary to that of Remark
%\ref{firstremark},
%that given
%$\bsb{f}\in \widetilde{\mathcal{M}}_{\mu,\nu,m}^{\leq l}$, we have $\chi^{k}\bsb{f}\in \widetilde{\mathcal{M}}[1/h]_{\mu,q^k\nu,m}^{\leq lq^k}$ for all $k\geq 0$.
%For example, we have seen that $\chi\bsb{E}=-\bsb{h}_q$.
%Formula (\ref{splitting}) relates in some sense the transcendence methods of Anderson-Brownawell-Papanikolas in \cite{ABP}
%and Mahler's method (see the forthcoming paper \cite{Pe2}). We believe that the study of the operator $\chi$ on almost-$A$-quasimodular forms is also worth to be pursued.}
%\end{Remarque}

\begin{Remarque}
{\em It is natural to ask whether any of the finitely generated $\TT_{>0}$-modules we have introduced so far, like $\mathcal{M}_{\mu,\nu,m},\MM_{\mu,\nu,m}$,
can be endowed with a natural extension of Hecke operators as defined, on the vector spaces $M_{w,m}$, in \cite{Ge}.
Let us fix $(\mu,\nu,m)$ in $\ZZ^2\times\ZZ/(q-1)\ZZ$.
One is then tempted to choose, for $\mathfrak{P}$ a prime ideal of $A$ generated by the monic polynomial
$\mathfrak{p}\in A\setminus\{0\}$ of degree $d>0$, a $\TT_{>0}$-linear map
$$T_{\mathfrak{P}}:\mathcal{M}_{\mu,\nu,m}\rightarrow\boldsymbol{\mathcal{O}}$$ by:
$$(T_{\mathfrak{P}}\bsb{f})(z):=\mathfrak{p}^\mu\overline{\mathfrak{p}}^\nu \bsb{f}(\mathfrak{p}z)+\sum_{b\in A,\deg_\theta b<d}\bsb{f}\left(\frac{z+b}{\mathfrak{p}}\right),$$
for $\bsb{f}\in\mathcal{M}_{\mu,\nu,m}$.

This map reduces to the Hecke operator $T_{\mathfrak{P}}:M_{w,m}\rightarrow M_{w,m}$ of \cite{Ge}
if we set $t=\theta$ (when this operation makes sense). However, it is unclear when the image $T_{\mathfrak{P}}(\mathcal{M}_{\mu,\nu,m})$ is
contained in $\mathcal{M}_{\mu,\nu,m}$. For instance, it can be proved that $\Delta^2$ is not eigenform for all the operators $T_{\mathfrak{P}}$.
Since it is easy, applying the second inclusion in (\ref{doubleinclusions}), to prove that if $q\neq 2,3$ then $\mathcal{M}_{2(q-1),2q(q-1),0}=\TT_{>0}\bsb{h}_q^{2(q-1)}$,
it follows that $T_{\mathfrak{P}}(\mathcal{M}_{2(q-1),2q(q-1),0})\subsetneq\mathcal{M}_{2(q-1),q(q-1),0}$.

The problem arises with our second factor of automorphy $(\bsb{J}_{\gamma})_{\gamma\in\Gamma}$:
it does not extend to a factor of automorphy $\mathbf{GL}_2(K)\times\Omega\rightarrow C[[t]]$. In other words, 
almost $A$-quasi-modular forms do not always come from lattice functions.
In this sense, 
our algebra $\mathcal{M}$ might still be ``too small", needing to be embedded in a larger algebra which is presently unknown.}
\end{Remarque}


\begin{thebibliography}{99}

\bibitem{anderson} G. Anderson, {\em $t$-motives.} Duke Math. J. 53, 457-502 (1986).

\bibitem{ABP} G. Anderson, D. Brownawell, M. Papanikolas, {\em Determination of 
the algebraic relations among special $\Gamma$-values in positive characteristic.} Ann. of Math. 160 
(2004), 237--313.

\bibitem{BP} V. Bosser \& F. Pellarin. {\em Differential properties of Drinfeld quasi-modular forms.}
Int. Math. Res. Not. 2008, No 11, Article ID rnn032, 56 p. (2008).

\bibitem{BP2} V. Bosser \& F. Pellarin. {\em On certain families of Drinfeld quasi-modular forms.}
J. Number Theory 129, No. 12, 2952--2990 (2009).

\bibitem{cor} G. Cornelissen. {\em A survey of Drinfeld modular forms.} Proceedings of the workshop on
Drinfeld modules, modular schemes and applications, Gekeler (ed.) et al., World Scientific, 167--187 (1997).

\bibitem{Go1} D. Goss. {\em $\pi$-adic Eisenstein series for Function Fields.} Compositio Math. 41, pp. 3--38 (1980).

\bibitem{Ge}  E.-U. Gekeler. {\em
On the coefficients of Drinfeld modular forms.}
Invent. Math. 93, No.3, 667--700 (1988).

\bibitem{gekeler:compositio} E.-U. Gekeler. {\em Quasi-periodic functions and Drinfeld
modular forms.} Compositio Math. 69, No. 3, 277--293 (1989).

\bibitem{Pa} M. A. Papanikolas. {\em Tannakian duality for Anderson-Drinfeld
motives and algebraic independence of Carlitz logarithms.},
Invent. Math. 171, 123--174 (2008).

\bibitem{Bourbaki} F. Pellarin. {\em Aspects de l'ind\'ependance alg\'ebrique en caract\'eristique non nulle.}
S\'eminaire Bourbaki Mars 2007 59me ann\'ee, 2006--2007, no. 973


\bibitem{Pe} F. Pellarin. {\em Estimating the order of vanishing at infinity of Drinfeld quasi-modular forms.} arXiv:0907.4507

\bibitem{Pe2} f. Pellarin. {\em An introduction to MahlerÕs method for 
transcendence and algebraic independence.} Manuscript, (2010)

\end{thebibliography}
\end{document}